\documentclass[aps,preprint,showpacs,amsmath,nofootinbib]{revtex4}
\usepackage{epsfig}
\usepackage{amssymb}
\usepackage{multirow}
\usepackage{dcolumn}
\usepackage{bm}
\usepackage[nottoc]{tocbibind}
\usepackage[ruled,vlined,onelanguage]{algorithm2e}

\newtheorem{proof}{Proof}
\newtheorem{theorem}{Theorem}
\newtheorem{proposition}{Proposition}
\usepackage{lineno,hyperref}
\modulolinenumbers[5]
\usepackage[utf8]{inputenc}
\usepackage{amsmath}
\modulolinenumbers[1]
\usepackage{siunitx}
\usepackage{amsmath}
\usepackage[capitalise]{cleveref}
\usepackage{siunitx}

\begin{document}

\title{Generalized Hopf Bifurcation in  a Cancer Model with Antigenicity under Weak and Strong Allee Effects}

\author{Eymard Hern\'andez-L\'opez$^1$}
\author{Mayra N\'u\~nez-L\'opez$^2$} 
\email{mayra.nunez@itam.mx}
\author{Napole\'on Navarro-Tito$^1$}

\affiliation{$^1$ Facultad de Ciencias Qu\'imico-Biol\'ogicas, Universidad Aut\'onoma de Guerrero, Chilpancingo, Guerrero 39000, M\'exico}

\affiliation{$^2$ Department of Mathematics, Instituto Tecnol\'ogico Aut\'onomo de M\'exico,
R\'io Hondo 1, Col. Progreso Tizap\'an, Ciudad de M\'exico,
C.P. 01080, M\'exico}

\begin{abstract}
This article deals with an autonomous differential equation model that studies the interaction between the immune system and the growth of tumor cells with strong and weak Allee effects. The Allee effect refers to interspecific competition, and when the population is small, it can retard population growth. The work focuses on describing analytically, using a set of parameters, the conditions in the phases of the immunoediting theory, particularly in the equilibrium phase, where a latent tumor would exist. Saddle-Node, Saddle-symmetric, Hopf, generalized Hopf, and Takens-Bogdanov bifurcations get presented for both Allee effects, and their biological interpretation regarding cancer dynamics gets discussed. The Hopf and generalized Hopf bifurcation curves get analyzed through hyper-parameter projections of the model, where it gets observed that with a strong Allee effect, more tumor control persists as it has higher antigenicity, in contrast to the weak Allee effect, where lower antigenicity gets observed. Also, we observe that the equilibrium phase persists as antigenicity increases with a strong Allee effect. Finally, the numerical continuation gets performed to replicate the analytical curves' bifurcations and draw the limit and double limit cycles.
\\
\\
{\it Keywords: Generalized Hopf bifurcation, Cancer modeling, Immunoediting, Weak-strong Allee effects.} 
\end{abstract}

\maketitle

\section{Introduction}
Multiple factors, such as genetic, environmental, viral infections, behavioral, and dietary factors, could cause the appearance of cancer. A latent tumor (equilibrium phase) is one of the crucial stages in cancer dynamics. Hence, the immune system is crucial in controlling and eliminating cancerous tumors. Cancer vulnerabilities become more apparent if thinking about the cells that make up a tumor as an endangered species. According to \cite{Altrock2015}, many tiny tumors possibly get formed but almost always become extinct before they are clinically relevant and even before possible detection. Extinction is a complex phenomenon often driven by the interaction of ecological and evolutionary processes; in this context, the Allee effect implies the existence of a growth threshold that may be explored in therapeutics \cite{Kolev}.

On the other hand, in the presence of Allee effects, successful tumors occur following rare large fluctuations in the population size that take the tumors over the Allee threshold. Allee effects get caused by several mechanisms, including cooperative feeding and defense, which can potentially be relevant to cancer. When the number of organisms is small, cooperation is typically inefficient, leading to a growth threshold in the population. Cell cooperation might be required in diseases like cancer to produce a sufficient density for tumor proliferation. Allee effect in cancer is being investigated because various factors may limit growth in tumor cells at low densities, cooperation among cells might be required to produce a sufficient density of diffusible growth factors needed for tumor proliferation \cite{Gore2009}, \cite{Axelrod2015}.

The Allee effect receives attention in several mathematical models, such as delay in differential equations where population density changes, in time series for population processes, and areas of biomathematics, among others \cite{Courchamp2008}. The weak Allee effect is when at low density, the per capita growth rate is lower than at high densities, but the growth rate remains positive even in this case of a small population. In contrast, the strong Allee effect is characterized by the fact that after going from a positive to negative per capita growth rate in populations, through the threshold curve, the per capita rate of population growth drops dramatically, becoming negative at an accelerated rate until the population's extinction. 

Simplified models at the cellular level consider tumor and immune cells in the microenvironment, where cytokines regulate their interactions and dormancy patterns \cite{Bellomo}, \cite{Kolev}. In literature, there are several dynamic models based on the interaction between cancer cells and the immune system, as in  \cite{Boer1985}, \cite{Kuznetsov1994}, \cite{Delgado2020}. 

In others, it's analized the interaction between the immune system’s cells and tumor cells has already been studied without Allee effect by
\cite{Kirschner1998}, \cite{Starkov2013},  \cite{Xavier2017}, \cite{Wei2013}. Most of the works mentioned above focus on stability
points and bifurcations found numerically, but
analytical bifurcations are lacking.
In \cite{hernandez2021bifurcations}, the authors modeled the interaction between tumor cells and the immune system’s effector cells; they study the generalized Hopf bifurcations (known as Bautin bifurcation), including weak Allee effect. Finally, in \cite{Rocha2019}, present another way to analyze Allee’s effects of bifurcation; they realized a study of the Allee effect on a generalized logistic map that exhibits rich and complex dynamics.

In this work, we propose a mathematical model to study of bifurcations on the dynamics which considers strong and weak Allee effects as a limiting factor in tumor growth and its interaction with the immune system.

The paper gets organized as follows: in Section \ref{sectionModel}, we present and describe the model with weak and strong Allee effect; in Section \ref{SectionCriticalPoint}, we present the critical points. Then, section \ref{SectionBifurcations} presents the bifurcations with the Allee effect. Finally, in Section \ref{SectionConclusion}, we draw some conclusions about this work.

\section{Model}\label{sectionModel}

The immune system is the first line of defense against aggressions of external nature, when normal cells become cancer cells,
the immune system launches. In this works we will study  the dynamics between cells of the immune system $E$ and tumor cells $T$
which will not take
into account the action of some treatment on cancer. In this system we consider weak and strong Allee effects in a single model so as to allow for a  comparative study of impacts of these types of Allee effect therefore the system is given by 

\begin{equation}\label{Panne_Jr}
\begin{array}{lcll}
      \frac{dT}{dt}&=&\frac{\hat{r} T (1-\hat{b} T)  (T-\hat{\beta} )}{\hat{\alpha} +T}-\frac{\hat{a} T E}{\hat{g}+T} &\\
    \frac{dE}{dt}&=&\hat{c} T-\hat{\mu}  E,
\end{array}
\end{equation}

\noindent where $\hat{c}$, $\hat{r}$, $\hat{a}$, $\hat{g}$, $\hat{\mu}$ $>0$.

In the model, the Allee threshold and the carrying capacity of the environment  explicitly included as model parameters.

The parameter $\hat{r}$ is the intrinsic growth rate of $T$, it can be understood as the proliferation rate
of tumor cells, and $1/\hat{b}$ is the carrying capacity of $T$. A logistic tumor cell growth represents the first
term of Eq. (1). The strong Allee threshold is given by  $\hat{\beta}$ and $\hat{\alpha}$ is a ``control" parameter that allows us to transition between weak effect Allee and strong effect Allee. If $\hat{\beta}=-\hat{\alpha}$ is the scenario in which the demographic Allee effect is absent. For our analysis, if $\hat{\beta}>0$ this corresponds to a strong effect Allee scenario, on the other hand, if $\hat{\beta}=0$, we obtain a weak effect Allee scenario with $\hat{\alpha}>0$ according to \cite{Capistran2018}, \cite{Courchamp2008}, \cite{Celik2009}.

The second term of first equation of system (1) represents the loss of
tumor cells due to their interaction with the effector cells $E$ at a rate of $\hat{a}$. The parameter $\hat{a}$ corresponds to the immune system’s response capacity
to the presence of tumor cells where $\hat{g}$ is the half-saturation rate for cancer clearance.

The second equation models the change in the population of the immune system’s effector cells with respect to time. The first term  represents the recruitment of effector cells in response to tumor antigenicity, where $\hat{c}$ is the tumor antigenicity (the
ability of the tumor to elicit an immune system response). Finally, the second term represents the death or apoptosis of the effector cells with $\hat{\mu}$  as death rate of immune cells.

To facilitate analysis,
we will use the following dimensionless system to perform the bifurcation analysis of the system (\ref{Panne_Jr}).
\begin{equation}\label{Panne_Jr_a}
\begin{array}{lcll}
      \frac{dx}{dt}&=&\frac{r x (1-b x)  (x-\beta )}{\alpha +x}-\frac{a x y}{g+x} &\\
    \frac{dy}{dt}&=&c x-\mu  y.
\end{array}
\end{equation}

\noindent where $x$ are the tumor cells and $y$ are cells of the
immune system.

\section{Critical Points}\label{SectionCriticalPoint}
The critical points for system (\ref{Panne_Jr_a}) with strong Allee effect, are the trivial point $P_{0}=(0,0)$, two complex critical points $P_{1},\bar{P_{1}}$, and real critical point given by 
\begin{equation}\label{PC1}
    P_{2}=\left(\frac{\eta -2 a c}{6 b \Gamma  \mu  r},\frac{c (\eta -2 a c \Gamma )}{6 b \Gamma  \mu ^2 r}\right),
\end{equation}

The equilibrium points of interest are those such that their inputs are real nonnegative. That is, the equilibrium point of interest is $P_2$,  $\Gamma$ and $\eta$ parameters can be found in the  Appendix \ref{Ape1} as (\ref{GAMMA})  and  (\ref{ETA}) respectively. 

\begin{proposition}
    All solutions $P_{2}=(x_0,y_0) \in \mathbb{R}^2_{+}$ from (\ref{PC1}) of system (\ref{Panne_Jr_a}) are positive in  the following strong Allee case.
        \begin{equation}\label{StrongA}
            \begin{array}{lcll}
              W_{Strong}&=&\left\{ W_{s_1}<\alpha<W_{s_2}, 0<a<W_{s_3}, W_{s_4}<\Gamma,0<b<W_{s_5} \right\}, 
            \end{array}
        \end{equation}
       where $0<\Gamma\leq 1$ with $\eta > 2 a c$, for $g=\beta$, and provided that all other parameters are positive, including the parameters $W_{s_{i}}>0$
       
\footnote{$W_{s_{1}}=\frac{a^2 c^2-2 a c \mu  r+b^2 \beta ^2 \mu ^2 r^2+\mu ^2 r^2}{a b c \mu  r}$, $W_{s_2}=\frac{\mu  r \left(\sqrt[3]{2} \mu  r \left(3 b^2 \beta ^2+1\right)+\Gamma \right)-a \left(2 \sqrt[3]{2} c \mu  r+c\right)}{3 \sqrt[3]{2} a b c \mu  r}$, $W_{s_{3}}=\frac{3 \sqrt[3]{2} b^2 \beta ^2 \mu ^2 r^2+\sqrt[3]{2} \mu ^2 r^2+\Gamma  \mu  r}{2 \sqrt[3]{2} c \mu  r+c}$, $W_{s_{4}}=\frac{1}{3} \sqrt{\frac{7 (a c-\mu  r)^3}{3 \left(\beta ^2 \mu ^3 r^3\right)}}$, and $W_{s_{5}}=\frac{3 \sqrt[3]{2} a^2 c^2-4 \sqrt[3]{2} a c \mu  r+a c+2 \sqrt[3]{2} \mu ^2 r^2}{\mu  r}$.}, for $i\in\{1,...,5\}$.
   
\end{proposition}

\begin{proof}
    If we consider the expression (\ref{GAMMA}), when simplifying it is enough to compare 

    \begin{eqnarray*}
        3 a^2 c^2 \mu  r (b (3 \alpha +2 \beta -2 g)+2)&>&2 a^3 c^3+3 a c \mu ^2 r^2 \left(b \left(3 \alpha +b \beta  (3 \alpha +2 \beta )+2 b g^2-g (b (3 \alpha +\beta )+1)+\beta \right)\right) \\
        &+&\mu ^3 r^3 (b (g-\beta )+2) (b (\beta +2 g)+1) (b (2 \beta +g)-1),
    \end{eqnarray*}

    If we simplify the expression for $\eta$ (\ref{ETA}), and the we have 
    $$\mu  r \left(\sqrt[3]{2} \mu  r \left(b \left(\beta  (b \beta -1)+b g^2+b \beta  g+g\right)+1\right)+b \Gamma  (\beta -g)+\Gamma \right)>a c \left(\sqrt[3]{2} \mu  r (3 \alpha  b+2 b \beta -2 b g+2)+1\right).$$

    The expression (\ref{PC1}) have positive components if 
    $$(0<\Gamma \leq 1\land a>0\land \eta >2 a c), \hspace{0.5cm}or \hspace{0.5cm} \Gamma >1\land a>0\land \eta >2 a c \Gamma.$$ 

    Therefore, the parameters $\Gamma$ and $\eta$ are positives if  we have the following conditions  $W_{s_1}<\alpha<W_{s_2}, 0<a<W_{s_3}, W_{s_4}<\Gamma,0<b<W_{s_5}$, for $\beta=g$, and the positive parameters of system (\ref{Panne_Jr_a}).
\end{proof}

\begin{figure}\centering
 \hspace{0.005cm} \textbf{(a)} \hspace{13.3cm} \textbf{(b)}
\includegraphics[scale=0.4]{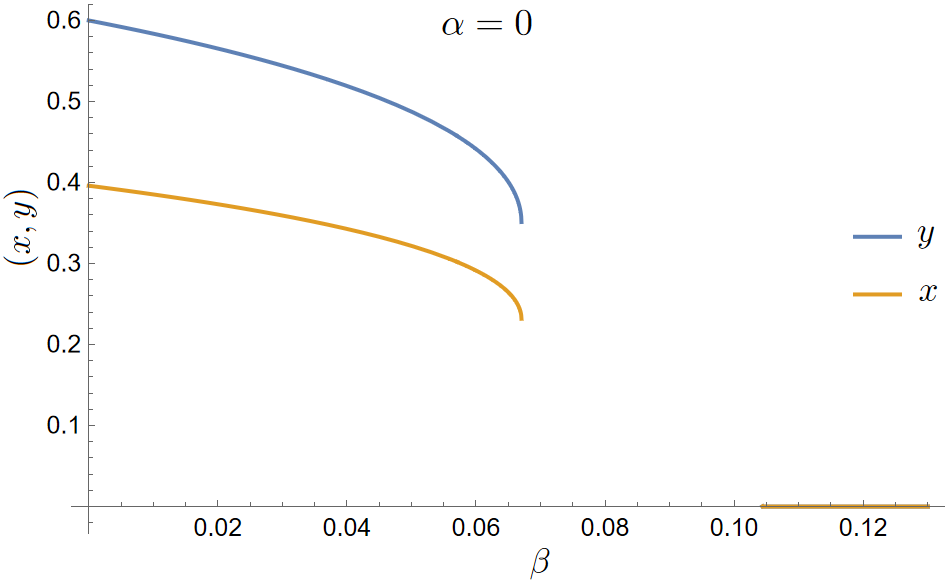}
\includegraphics[scale=0.4]{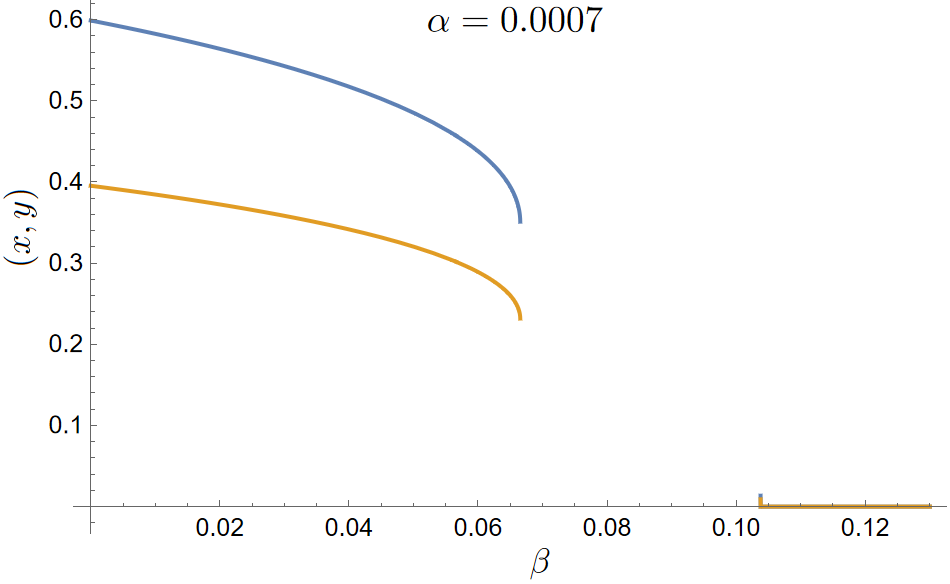}
\\
 \hspace{0.005cm} \textbf{(c)} \hspace{13.3cm} \textbf{(d)}
\\
\includegraphics[scale=0.4]{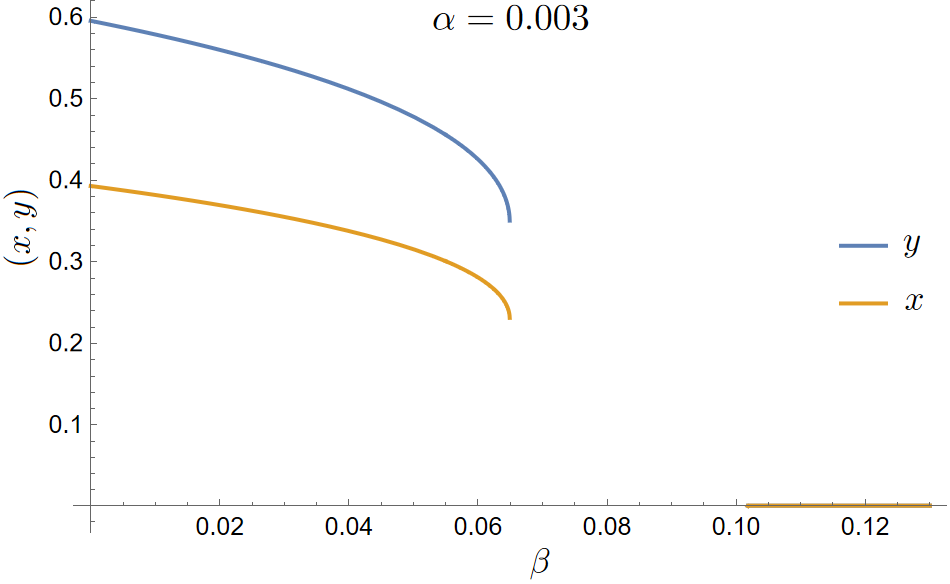}
\includegraphics[scale=0.4]{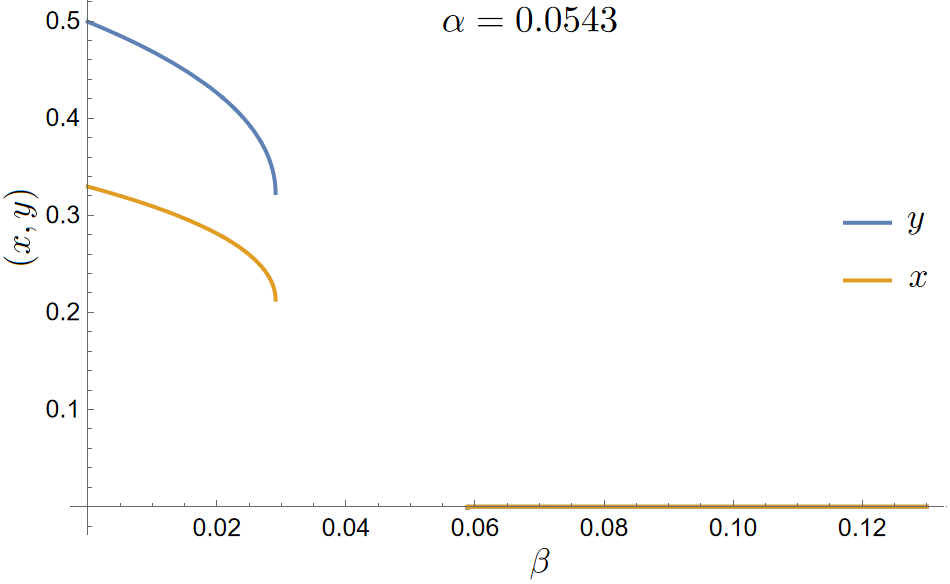}
\caption{The effect of parameters $\alpha$ and $\beta$ on the values of the critical point (\ref{PC1}).}
\label{fig:PunCrit1}
\end{figure}

\begin{figure}[htb!]
\begin{center}
\includegraphics[scale=0.5]{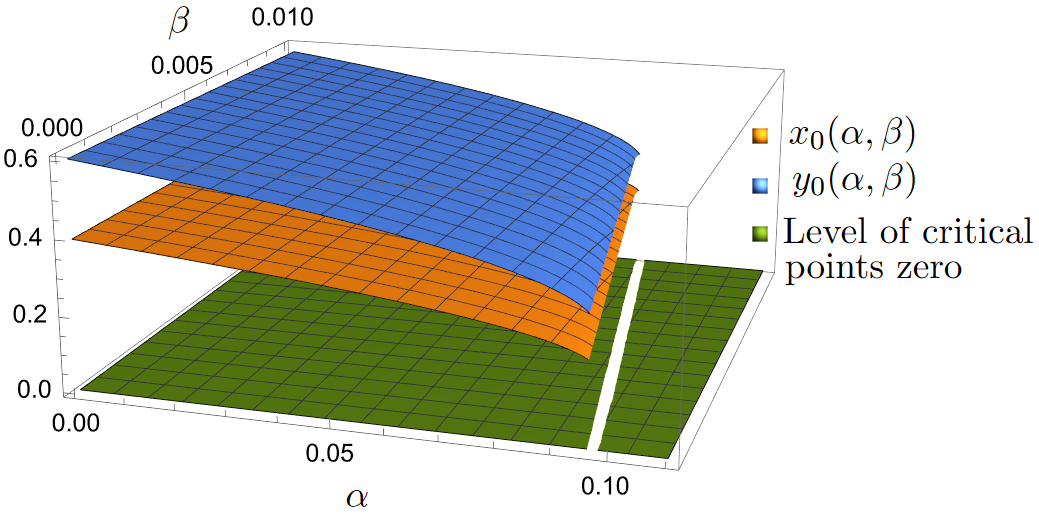}
\end{center}
\caption{The effect of parameters $\alpha$ and $\beta$ on the values of the critical point (\ref{PC1}), for parameter values $a=5.55556$, $b=1.24997$, $c=0.252521$, $g=0.0001$, $r=16.6667$, and $\mu = 0.166667$.}
\label{fig:PunCrit2}
\end{figure}

 Figure \ref{fig:PunCrit1} shows the behavior of critical point $P_2$ given by expression (\ref{PC1}) with admissible parameter values, the parameters corresponding to
 different parameterizations of strong Allee effect, represented by $\alpha$ and $\beta$, 
 we show the critical values of $x$ and $y$ for different values of $\alpha$.
As $\alpha$ increases  the population density at which the per capita growth rate is maximized moves to a higher density, $i.e.$ reaches lower maximum values.
Also for the case of Figure 1 (b)-(d), when $\beta=0$ a weak Allee effect is obtained, just the maximum point of the curves. Figure \ref{fig:PunCrit2} shows curves in 3D representing the critical points by varying $\alpha$ and $\beta$. If the value of $\alpha$ increases, the parameter $\beta$ should be smaller.
\\
\\

We analyze the stability of the critical points of biology interest in the following result.
\begin{proposition}\label{Pcritical}
If the critical points are reals, we have the stability as follows.
\begin{description}
\item[(a)] The critical point $P_{0}$ is asymptotically stable.
\item[(b)] The critical point $P_{2}$ is an unstable or asymptotically stable node, unstable or asymptotically stable spiral, stable center, or a saddle.
\end{description}
\end{proposition}
\begin{proof}
Analyzing the eigenvalues for $P_{0}$: $E_{0_1}=- \beta  g r$, and $E_{0_2}=-\alpha  g \mu$, then the point $P_{0}$ is asymptotically stable. We consider the eigenvalues for $P_{2}$ as follow 
$$E_{2_1,2_2}=\frac{1}{2}\left(\hat{P_{2_1}}-\hat{P_{2_2}}\right)\pm\sqrt{F(x_0,y_0:a,b,c,g,r,\alpha,\beta)}),$$
be 
$$
\Omega=(x_0,y_0:a,b,c,g,r,\alpha,\beta)
$$
such that
$F(\Omega)=F(x_0,y_0:a,b,c,g,r,\alpha,\beta)$, and $\hat{P_{2_1}}=(r x_0 (b \beta +1) (2 g+3 x_0)$, $\hat{P_{2_2}}=\left(a y_0 (\alpha +2 x_0)+b r x_0^2 (3 g+4 x_0)+\beta  r (g+2 x_0)+\mu  (g+x_0) (\alpha +x_0)\right)$, then
\begin{itemize}
    \item For $F(\Omega)=0$: If $\hat{P_{2_1}}>\hat{P_{2_2}}$, $P_{2}$ is an unstable node. But if $\hat{P_{2_1}}<\hat{P_{2_2}},$ then $P_{2}$ is asymptotically stable.

    \item In case $F(\Omega)<0$: if $\hat{P_{2_1}}>\hat{P_{2_2}}$, $P_{2}$ is spiral unstable. But if $\hat{P_{2_1}}<\hat{P_{2_2}},$ $P_{2}$ is spiral asymptotically stable, either if $\hat{P_{2_1}}=\hat{P_{2_2}}$ $P_{2}$ is stable center.

    \item For $F(\Omega)>0$, $P_{2}$ is an unstable node if $\hat{P_{2_1}}>\hat{P_{2_2}}$ and $(\hat{P_{2_1}}-\hat{P_{2_2}})\pm \sqrt{F(\Omega)}$, or if $\hat{P_{2_1}}<\hat{P_{2_2}}$ and $(\hat{P_{2_1}}-\hat{P_{2_2}})\pm \sqrt{F(\Omega)}$. The point $P_{2}$ is a saddle if $\hat{P_{2_1}}>\hat{P_{2_2}}$, $(\hat{P_{2_1}}-\hat{P_{2_2}})- \sqrt{F(\Omega)}<0$, and  $(\hat{P_{2_1}}-\hat{P_{2_2}})+ \sqrt{F(\Omega)}>0$, or if $\hat{P_{2_1}}<\hat{P_{2_2}}$, $(\hat{P_{2_1}}-\hat{P_{2_2}})- \sqrt{F(\Omega)}<0$, and  $(\hat{P_{2_1}}-\hat{P_{2_2}})+ \sqrt{F(\Omega)}>0$. Finally, $P_{2}$ is asymptotically stable if $\hat{P_{2_1}}<\hat{P_{2_2}}$, and $(\hat{P_{2_1}}-\hat{P_{2_2}})\pm \sqrt{F(\Omega)}<0$.
    
\end{itemize}
The function $F(\Omega)$ in variables with parameters is in the Appendix \ref{Ape1} as (\ref{FPC1}).
\end{proof}

\begin{figure}[htb!]
\begin{center}
\includegraphics[scale=0.65]{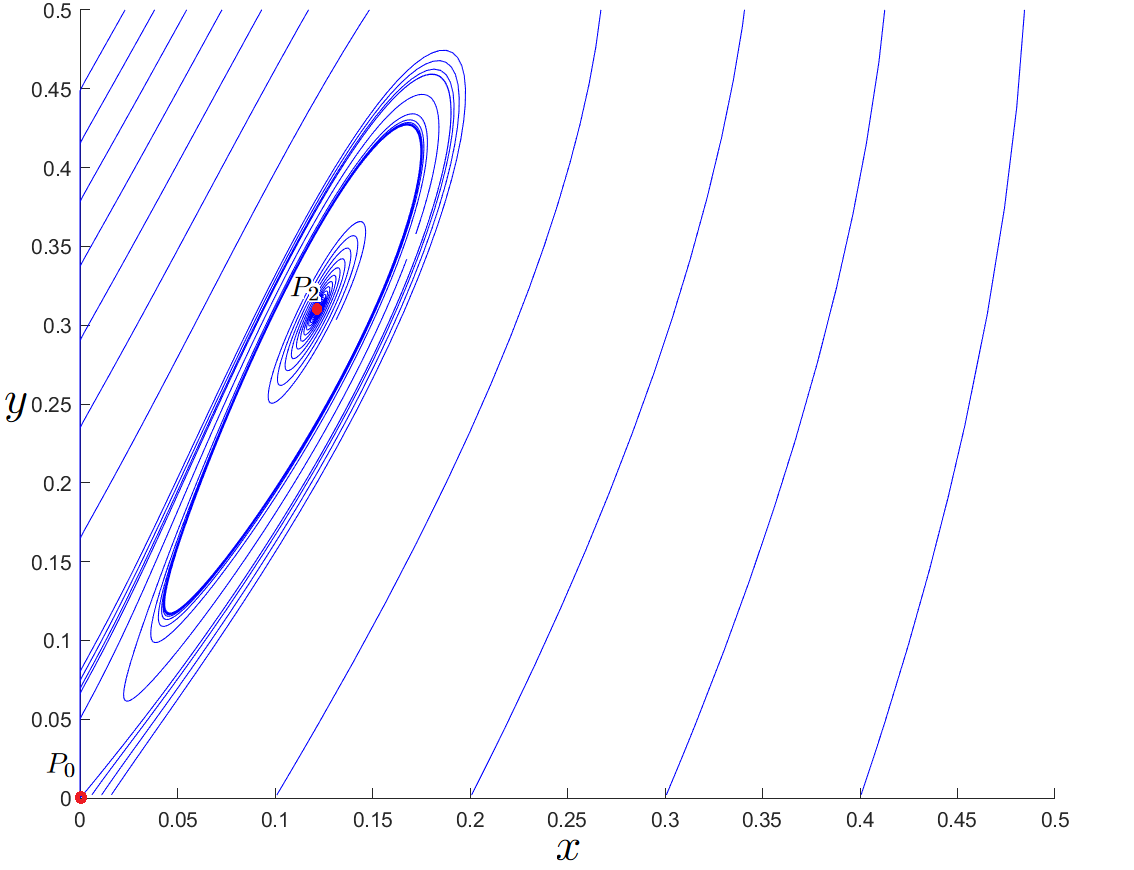}
\end{center}
\caption{Critical points and limit cycle for parameter values $a=5.55556$, $b=1.24997$, $c=30.4117$, $g=0.0001$, $r=16.6667$, $\alpha=1.0 \times 10^{-11}$, $\mu = 11.95$, and $\beta=1.0 \times 10^{-9}$.}
\label{fig1}
\end{figure}

Figure \ref{fig1}, shows the critical points $P_{0}$ and $P_{2}$, with some trajectories in the phase plane. It can be seen that there is a limit cycle for the set of positive parameters $a=5.55556$, $b=1.24997$, $c=30.4117$, $g=0.0001$, $r=16.6667$, $\alpha=1.0 \times 10^{-11}$, $\mu = 11.95$, and $\beta=1.0 \times 10^{-9}$, where  strong Allee effect gets represented by $\alpha$ and $\beta$.

\begin{table}[ht]
\caption {Parameter dimensional values of the System (\ref{Panne_Jr}).}
{\begin{tabular}{l c c}\\[-2pt]
\toprule
\textbf{Description}  &\textbf{Initials values} \\[6pt]
\hline\\[-2pt]
\textbf{$a$ : cancer clearance term limited} &$5.55556$\\[2pt]
\textbf{$1/b$ : carrying capacity of tumor cells}  & $1 \times 10^9$\\[2pt]
\textbf{$c$ : antigenicity}  & $30.4117$\\[2pt]
\textbf{$\mu$ : death rate of immune cells} & $11.95$\\[2pt]
\textbf{$r$ : cancer growth rate} & $16.6667$ \\[2pt]
\textbf{$g$ : half-sat. for cancer clearance} & $0.0001$ \\[2pt]
\textbf{$\alpha$ : weak Allee effect constant} & $1\times10^{-11}$ \\[2pt]
\textbf{$\beta$ : strong Allee effect constant}  & $1\times10^{-9}$ \\[2pt]
{} & &\\[1pt]
\botrule
\end{tabular}}
\end{table}

\section{Bifurcations}\label{SectionBifurcations}
In this section, we introduce some crucial results on bifurcations present in the system (\ref{Panne_Jr_a}). We first present bifurcations with oncological interpretation. Such is the case of the Hopf bifurcation since it is a helpful tool to understand the occurrence of oscillatory behavior in the system (\ref{Panne_Jr_a}) and can help explain how complex patterns of behavior can arise due to changes in the system parameters. The conditions for a generalized Hopf or Bautin bifurcation get then stated; this bifurcation occurs when a stable limit cycle loses its stability as a parameter is varied.

We mention bifurcations in the system (\ref{Panne_Jr_a}) that, despite not having ontological interpretation because the parameters are not all positive, are vital to study the Hopf and Bautin bifurcations. These are the saddle-node or Takens-Bogdanov bifurcations.

\begin{proposition}\label{propositionSN1}
    The following set of parameters 
    \begin{equation}\label{S_N}
        \varsigma=\left\{(a, \alpha,  b, \beta,  c, g, \mu,  r)\mid SN\right\},
    \end{equation}
  contains the saddle-node bifurcations of the system (\ref{Panne_Jr_a}), and the expression for $SN$ is in the Appendix \ref{Ape1}.
\end{proposition}
\begin{proof}
Calculating the roots in common $R_1$, of $\frac{dx}{dt}=f_{1}(x,y)$ and $\frac{dy}{dt}=f_{2}(x,y)$, with respect to the variable $y$, and taking the discriminant of $R_1$, now with respect to the variable $x$, then we obtain the expression for $SN$ in (\ref{S_N_1}). 
\end{proof}

\begin{theorem}\label{HopfTeo}
    The following set of parameters 
    \begin{equation}\label{H}
        H=\left\{(a, \alpha,  b, \beta,  c, g, \mu,  r)\mid HOPF\right\},
    \end{equation}
  contains the symmetric-saddle, and Hopf bifurcations  of the system (\ref{Panne_Jr_a}).
\end{theorem}
\begin{proof}
Let $f_1(x,y)$ and $f_2(x,y)$, as before, we consider the Jacobian matrix $A$ of the system (\ref{Panne_Jr_a}), with the trace $trA$. Let $R_2$ be the roots in common between the trace $trA$ and $f_1(x,y)$, and $R_3$ be the roots in common between the trace $trA$ and $f_2(x,y)$ regarding the $y$ variable. Finally, the $HOPF$ expression gets found by calculating the common roots between $R_2$ and $R_3$, with respect to the variable $x$, so the HOPF expression depends only on parameters. We obtain the expression $HOPF$ in (\ref{H_1}). 
\end{proof}

\begin{theorem}\label{BTToe}
    The following set of parameters 
    \begin{equation}\label{BT}
        B=\left\{(a, \alpha,  b, \beta,  c, g, \mu,  r)\mid BT=(BT_1)(BT_2)(BTI) \right\},
    \end{equation}
  contains the Bogdanov-Takens bifurcations  of the system (\ref{Panne_Jr_a}).
\end{theorem}
\begin{proof}
Calculating the roots in common between the $SN$ and $HOPF$ expressions results in a polynomial $BT$ with more than thirty-five thousand components dependent on parameters and not variables, where $$BT_1=\left(a \alpha  c-\alpha^2 b \mu  r-\alpha  b \beta  \mu  r-\alpha  \mu  r-\beta  \mu  r\right)^2,$$
the $BT_2$ expression is in Appendix \ref{Ape1} as (\ref{BT_1}), and the $BTI$ have many terms. 
\end{proof}

 As mentioned above, the saddle-node and Taken Bogdanov bifurcations only make sense for positive parameters; however, in the $SN$ and $BT$ expressions, some parameters must be negative for these bifurcations to exist, then these cases lack oncology interpretations. Furthermore, using symbolic software, we can only handle two principal terms on $BT$; the intractable term $BTI$ contains more than eighteen thousand terms. The main result of this work gets stated in the following theorem.

\begin{theorem}\label{BAUTT}
    The set 
    \begin{equation}\label{BAUE}
        Bau=\left\{(a, \alpha,  b, \beta,  c, g, \mu,  r)\mid BAUT \right\},
    \end{equation}
     contains the Generalized Hopf bifurcation points.
\end{theorem}
\begin{proof}
    Using $(8.22)$ and $(8.23)$ from \cite{Kuznetsov1994}, p. $310$, we calculate the first and second Lyapunov. To find the expression that determines the generalized Hopf bifurcation, we calculate the roots in common between the first Lyapunov coefficient $l_1(0)$ and the HOPF set, thus obtaining the BAUT expression analytically, where BAUT is a polynomial in parameters with more than one hundred forty-one thousand terms.
     $$BAUT=\left(28000 a^{14} \alpha ^{10} b c^{14} g^3 r^3+ \cdot \cdot \cdot + 16 \alpha ^6 b^{12} \beta ^{10} g^8 \mu ^{19} r^{12}\right). $$ 
     
     The second Lyapunov coefficient gets calculated analytically, but due to its size (more than eighty thousand terms, but now containing the variable $x_0$), obtaining an analytical expression to determine the sign was impossible. In this case, we used numerical values of the parameters and variables in a neighborhood where we observed a generalized Hopf bifurcation in numerical continuation. Then numerically, the second Lyapunov coefficient was positive; even when the variable $x_0$ is zero, the second Lyapunov coefficient has the following positive expression.
    \begin{eqnarray*}
        l_2(0)_{x_0 \sim 0} &=&
a^2 \alpha ^{14} \beta ^3 c^2 g^{13} r^3 \left(a \alpha ^2 c+\beta  r^2 (\alpha +\beta ) (\alpha  b+1)\right)^2\\
&\times&\left(10 a \alpha ^2 c+3 \beta  r^2 (\alpha +\beta ) (\alpha  b+1)\right)    
    \end{eqnarray*}
    The above expression is positive because the parameters are positive.
\end{proof}

\begin{figure}[htb]
\begin{center}
\includegraphics[scale=0.7]{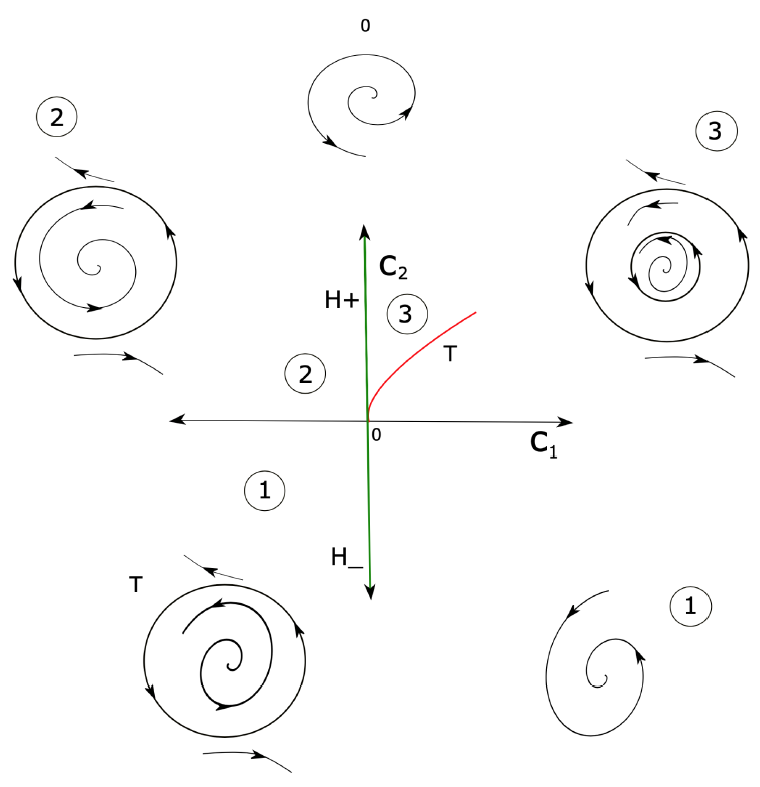}
\caption{Local Generalized Hopf bifurcation diagram of the system (\ref{Panne_Jr_a}). The bifurcation separates two branches, $H_{-}$ and $H_{+}$, corresponding to the Hopf bifurcation with a negative and positive first Lyapunov coefficient.  The parameters $C_1$ and $C_2$ indicate the Bautin bifurcation parameters in the model (\ref{Panne_Jr_a}), which are $c$ and $\mu$ respectively.}
\end{center}
\label{fig2}
\end{figure}

 Figure \ref{fig2} presents the local bifurcation diagram around the generalized Hopf point stated in Theorem \ref{BAUTT}; the diagram is similar to the one studied for the case without  Allee effect and with weak Allee  effect in \cite{hernandez2021bifurcations}. Although the first and second Lyapunov coefficients have extensive expressions, it was possible to study the main features and their sign in the case of the second coefficient. The existence of a generalized Hopf bifurcation can exhibit the persistence of double-limit cycles in a neighborhood. These double cycles characterize the equilibrium phase in cancer (latent tumor).

\begin{figure}[htb!]
\begin{center}
\includegraphics[scale=0.5]{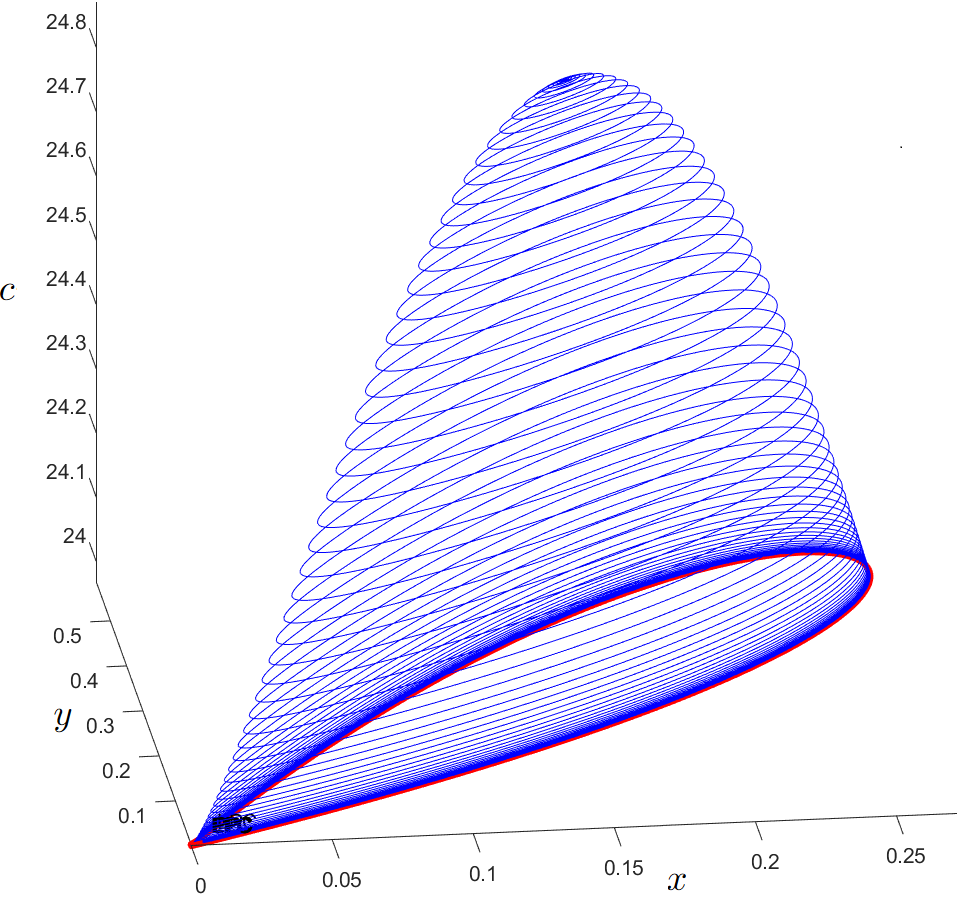}
\end{center}
\caption{Double limit cycle for parameter values $a=5.55556$, $b=1.24997$, $g=0.0001$, $r=16.6667$, $\alpha=1.0 \times 10^{-11}$, $\mu = 10.0530$, and $\beta=0$.}
\label{fig31}
\end{figure}

Figure \ref{fig31} shows the weak effect Allee scenario, this effect gets recovered when $\beta=0$ in Appendix Eq. (\ref{Bau_1}). The case when there is no Allee effect, neither strong nor weak, was studied in \cite{hernandez2021bifurcations}. Around the base of this graph, it can get seen that there are double limit cycles on values of the parameter $c$. The outer red cycle is unstable, and the inner blue cycle is stable, as shown by region 3 of the bifurcation diagram in Figure \ref{fig2}.
The occurrence of double cycles is presented for values very close to the specific value of $c = 24$,
marked by two red curves at the bottom of all the blue cycles  marked by two red curves at the bottom of all the blue cycles.

\begin{figure}[htb!]
\begin{center}
\includegraphics[scale=0.5]{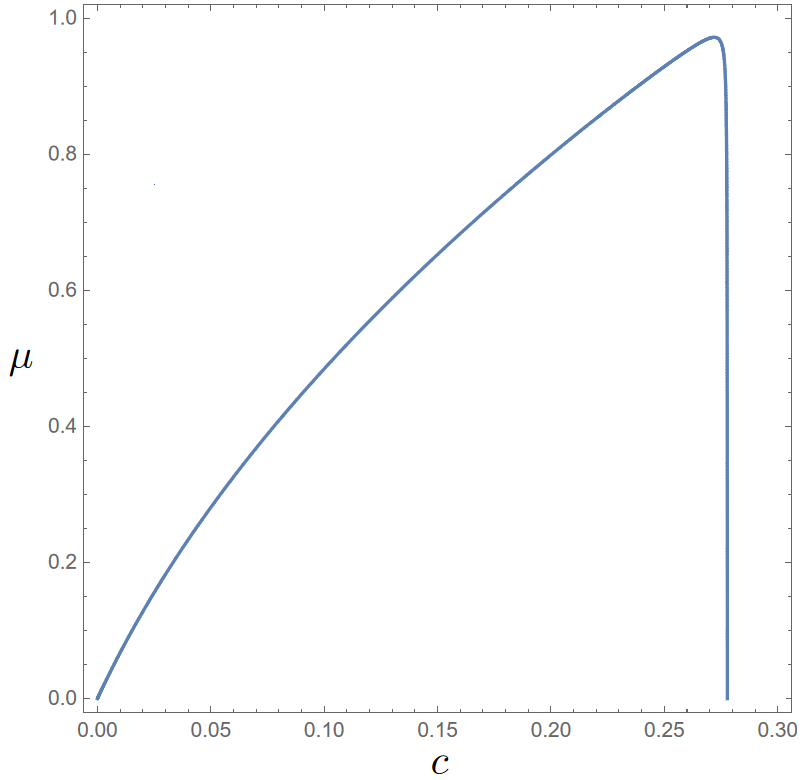}
\end{center}
\caption{Projection analytical of the level set of the symmetric-saddle and Hopf bifurcations (without Allee effect) from Theorem $1$ in \cite{hernandez2021bifurcations}, and here we recover it through our Theorem \ref{HopfTeo}, with parameters $a=3.59937$, $b=1$, $g=0.0001$, $r=1$, $\alpha=0$, and $\beta=0$.}
\label{fig4}
\end{figure}

\begin{figure}[htb!]
\begin{center}
\includegraphics[scale=0.5]{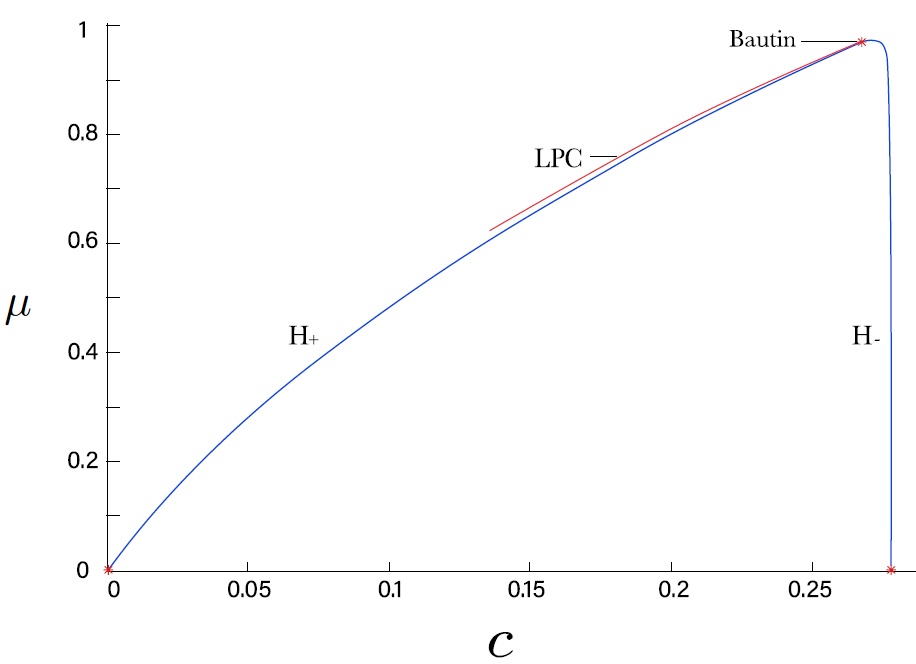}
\end{center}
\caption{Numerical continuation of symmetric-saddle and Hopf bifurcations (without Allee effect) with parameters $a=3.59937$, $b=1$, $g=0.0001$, $r=1$, $\alpha=0$, and $\beta=0$. The Hopf curve in blue color and the cycle limit point curve in red color.}
\label{fig41}
\end{figure}

\begin{figure}[htb!]
\begin{center}
\includegraphics[scale=0.5]{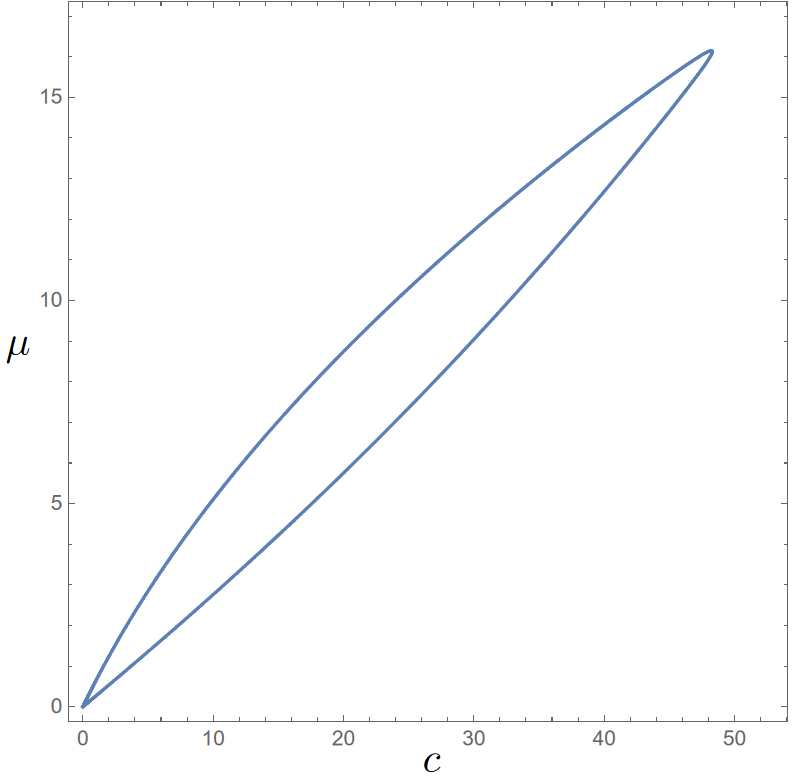}
\end{center}
\caption{Projection analytical of the level set of the symmetric-saddle and Hopf bifurcations (weak Allee effect) from Theorem $3$ in \cite{hernandez2021bifurcations}, and here we recover it through our Theorem \ref{HopfTeo}, with parameters $a=5.55556$, $b=1.24997$, $g=0.0001$, $r=16.6667$,  $\alpha=8.00021\times 10^{-5}$, and $\beta=0$.}
\label{fig5}
\end{figure}

\begin{figure}[htb!]
\begin{center}
\includegraphics[scale=0.7]{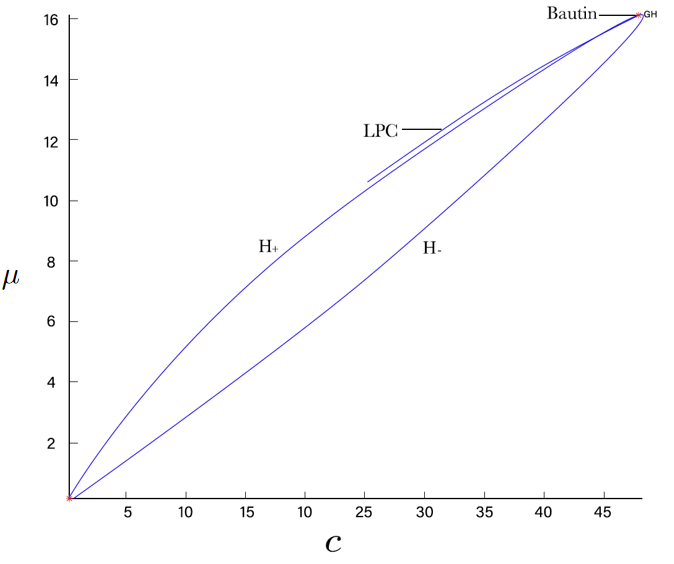}
\end{center}
\caption{Numerical continuation of symmetric-saddle and Hopf curve bifurcations (weak Allee effect). Hopf curve, LPC curve, and Bautin point. Parameters $a=5.55556$, $b=1.24997$, $g=0.0001$, $r=16.6667$,  $\alpha=1.0 \times 10^{-11}$, and $\beta=0$.}
\label{fig51}
\end{figure}

Figure \ref{fig4} shows a part of the level set corresponding to the saddle-symmetric and Hopf bifurcations analytically, as a function of the parameters $c$ and $m$, with a projection on the parameters values $a=3.59937$, $b=1$, $g=0.0001$, $r=1$, $\alpha=0$, and $\beta=0$. On the other hand, Figure \ref{fig41} is its counterpart in numerical continuation, and this figure corresponds to Figure $3$ of \cite{hernandez2021bifurcations}. Similar to the observations made previously, Figures \ref{fig5} and \ref{fig51} represent the case of the weak Allee effect; the first figure mentioned corresponds to the level set of saddle-symmetric and Hopf bifurcations analytically, while the second represents its numerical continuation, and in the same way, as in the previous case.

\begin{figure}[htb!]
\begin{center}
\includegraphics[scale=0.5]{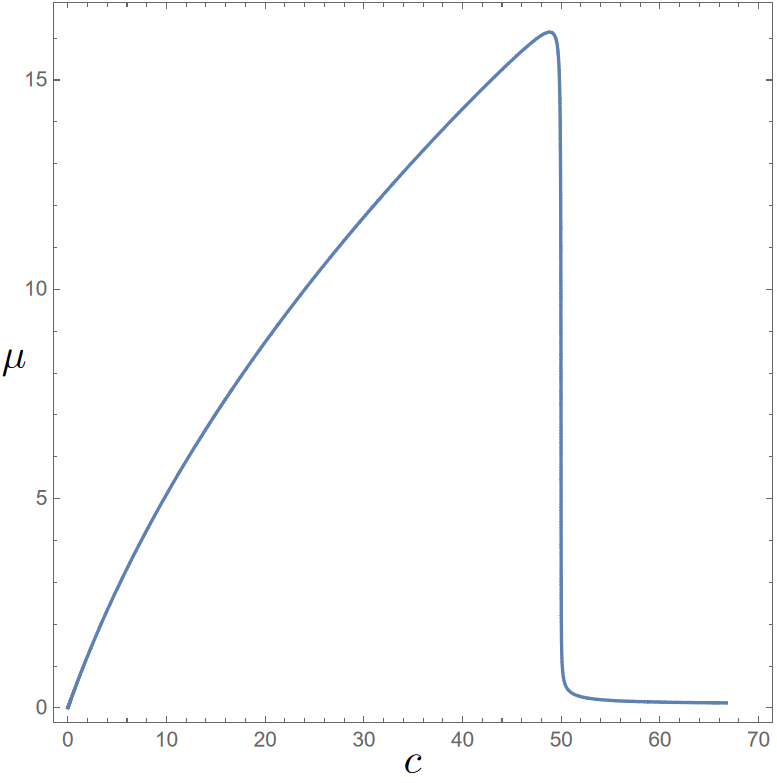}
\end{center}
\caption{Projection analytical of the level set of the symmetric-saddle and Hopf bifurcations (strong  Allee effect) from Theorem \ref{HopfTeo} with parameters $a=5.55556$, $b=1.24997$, $g=0.0001$, $r=16.6667$, $\alpha=1.0 \times 10^{-11}$, and $\beta=1.0 \times 10^{-9}$}
\label{fig6}
\end{figure}

\begin{figure}[htb!]
\begin{center}
\includegraphics[scale=0.5]{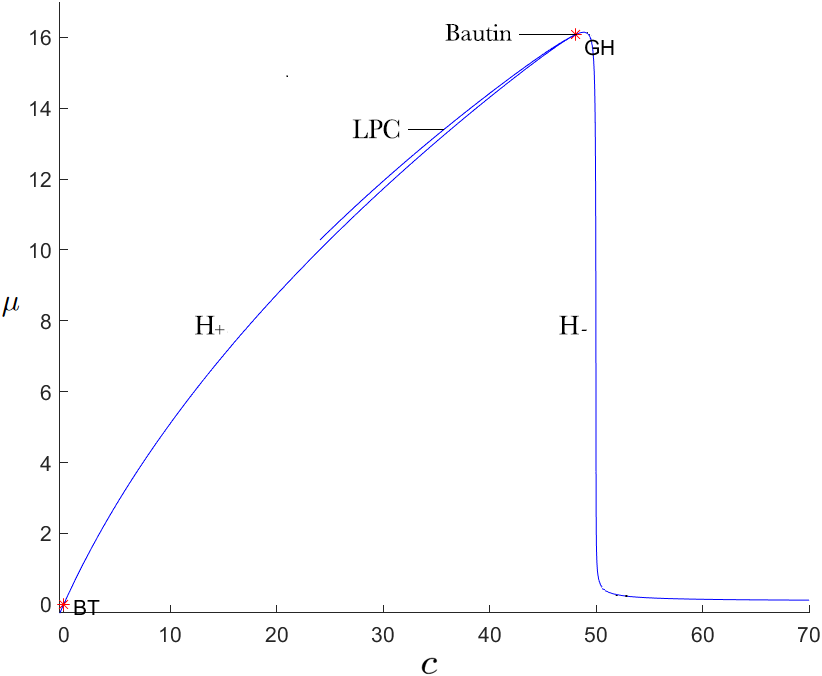}
\end{center}
\caption{Numerical continuation of symmetric-saddle and Hopf curve bifurcations (strong  Allee effect) with parameters $a=5.55556$, $b=1.24997$, $g=0.0001$, $r=16.6667$, $\alpha=1.0 \times 10^{-11}$, and $\beta=1.0 \times 10^{-9}$.}
\label{fig61}
\end{figure}

\begin{figure}[htb!]
\begin{center}
\includegraphics[scale=0.5]{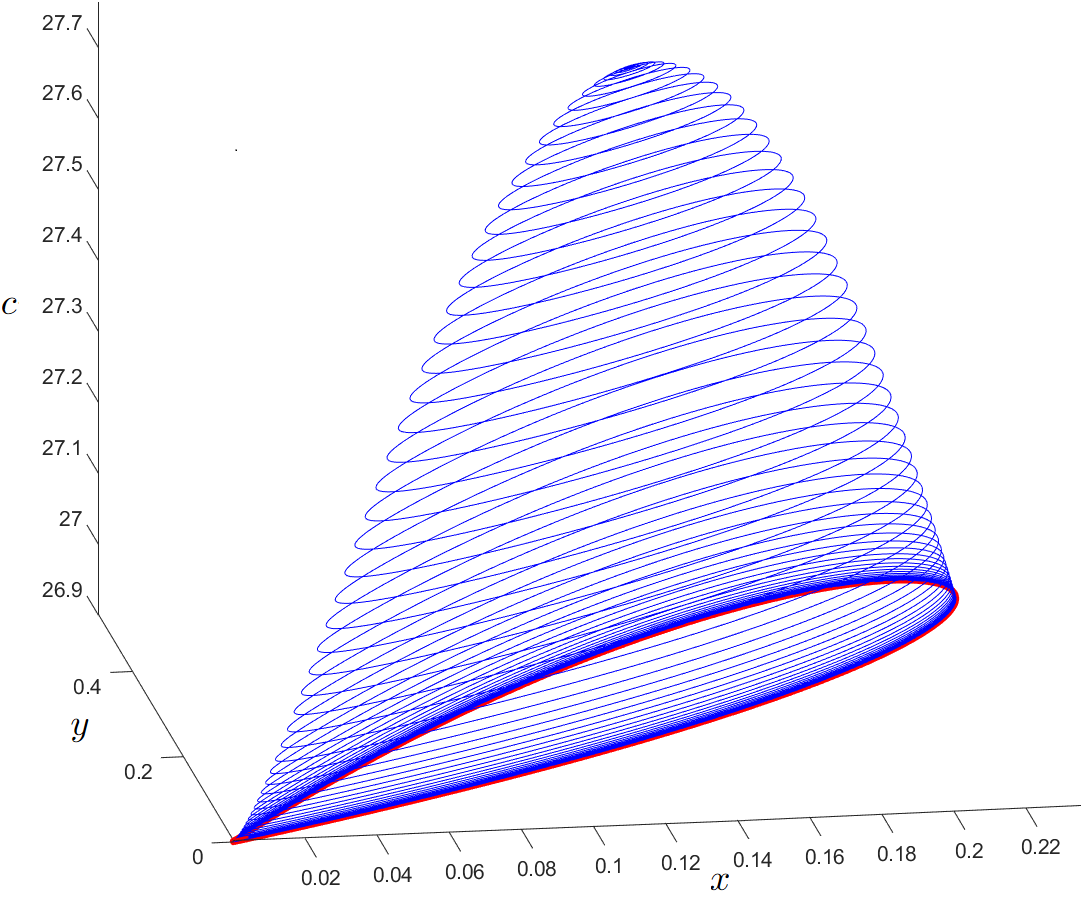}
\end{center}
\caption{Double limit cycle for parameter values $a=5.55556$, $b=1.24997$, $g=0.0001$, $r=16.6667$, $\alpha=1.0 \times 10^{-11}$, $\mu = 10.0530$, and $\beta=1.0 \times 10^{-7}$.}
\label{fig3}
\end{figure}

Part of the level set of Theorem \ref{BAUTT} is observed analytically in Figure \ref{fig6} for a strong Allee effect. The full expression is determined by the multiparametric polynomial (\ref{Bau_1}) and is plotted by a projection onto the parametric values $a=5.55556$, $b=1.24997$, $g=0.0001$, $r=16.6667$, $\alpha=1.0 \times 10^{-11}$, and $\beta=1.0 \times 10^{-9}$. The numerical continuation is the reflection of the analytical calculations. It gets shown in Figure \ref{fig61}, where the curve is a function of the parameters $c$ and $\mu$, the limit point of the cycles curve (LPC) get also observed, and the upper part of the Bautin bifurcation point get found using Matcont, as well as the LPC curve. 

The Bautin bifurcation we found we get characterized by the sign of the second Lyapunov coefficient and is a type of bifurcation where from an unstable limit cycle (Region $2$ in Figure \ref{fig2}), a new stable cycle emerges as the bifurcation parameters are perturbed, crossing the vertical  (Axis $C_2$ region $3$ in Figure \ref{fig2}), which subsequently collapses into the $LPC$ curve ($T$ curve in Figure \ref{fig2}), and a type of unstable limit cycle persists.

Figure \ref{fig3} shows the numerical continuation of double limit cycles  with a strong Allee effect; in this case, the occurrence of double cycles is presented for values very close to the specific value of $c=26.9$. By varying the parameter $c$, we can observe the presence of limit cycles with different periods with ranges of values for antigenicity in the interval $[26.9, 27.7]$ with $\beta=1.0 \times 10^{-7}$.

In Figures \ref{fig31} and \ref{fig3},  limit cycles get shown for different values of the parameter $c$; it can get seen that in the inner region of the paraboloid-shaped graph, there are limit cycles with a short period, which increases as the $c$ value decreases until reaching the cycle in red, which indicates a cycle limit point.  In addition, the weak and strong Allee effects get observed the difference in the growth intervals of tumor cells and the immune system cells since all the parameters are the same except for the value of $\beta$.

\section{Conclusion}\label{SectionConclusion}
This work analyzes a simple mathematical model to study the interaction between the immune system and tumor cell growth with strong and weak Allee effects. It is known that the Allee effect can represent interspecific competition; when the population is small, this is often reflected in a delay in population growth. An example of this is shown in Figure 1, where the tumor cell population’s growth rate decreases as parameters $\alpha$ and $\beta$ become more considerable, even becoming extinct. This work focuses on one of the phases of the immunoediting theory; this crucial stage in cancer is the equilibrium phase; in this phase, it gets considered that there is a latent tumor. Second, this work presents the saddle-symmetric and Hopf bifurcations analytically for both Allee effects, strong and weak, that correspond to the numerical continuations where the antigenicity of the tumor plays a key role.

Saddle-node bifurcation and Takens–Bogdanov bifurcation do not have a biological interpretation; however, the saddle-node curve helped calculate the Hopf bifurcation. The Hopf bifurcation indicates the existence of limit cycles. Bautin suggests two limit cycles; consequently, the model presents the equilibrium phase in immunoediting theory. In terms of cancer dynamics, a Bautin bifurcation determines another type of equilibrium, where double limit cycles occur, and therefore, it can also get characterized as the equilibrium phase. In this phase, a tumor keeps dormant, where tumor cells are unstable with mutations and resists the immune system (tumor cells use immune evasive strategies to grow and can be detected clinically, but in some cases, this can correspond to the escape phase). We can see that considering the same values of the system parameters, with a strong  Allee effect, there is greater tumor control by having a higher antigenicity, in difference with the weak Allee effect, where we can observe a smaller antigenicity (see Figures \ref{fig31} and \ref{fig3}).

According to the numerical continuation with a strong Allee effect, the equilibrium phase persists as the $c$ antigenicity increases (Hopf limit cycles). This is evidenced when the mortality rate $\mu$ tends to zero. However, compared with the numerical continuation without the Allee effect, the antigenicity no longer increases even though $\mu$ tends to zero.

\section*{ACKNOWLEDGMENTS} \noindent MN-L acknowledges the financial support from the Asociaci\'on Mexicana de Cultura, A.C. EH-L acknowledges the most valuable support from CONACYT (Consejo Nacional de Ciencia y Tecnolog\'ia) for a postdoctoral fellowship.

\newpage
\appendix
\section{}
\label{Ape1}
In this appendix we present several equations including the result of the calculations in Propositions and Theorems.  The following equation get obtained from the demonstration of Proposition \ref{Pcritical}. 

\begin{equation}\label{FPC1}
\begin{array}{lcll}
    F(x_0,y_0:a,b,c,g,r,\mu,\alpha,\beta)&=&16 b^2 r^2 x_0^6-24 b r^2 x_0^5+24 b^2 g r^2 x_0^5-24 b^2 r^2 \beta  x_0^5&\\
    &+&9 r^2 x_0^4+9 b^2 r^2 \beta ^2 x_0^4+\mu ^2 x_0^4-12 a c x_0^4+16 a b r y_0 x_0^4+34 b r^2 \beta  x_0^4&\\
    &+&6 r \mu  x_0^4-8 b r \alpha  \mu  x_0^4+6 b r \beta  \mu  x_0^4-12 b g^2 r^2 x_0^3+12 g r^2 x_0^3-12 b r^2 \beta ^2 x_0^3
    &\\
    &+&12 b^2 g r^2 \beta ^2 x_0^3+2 g \mu ^2 x_0^3+2 \alpha  \mu ^2 x_0^3-8 a c g x_0^3-12 a r y_0 x_0^3+12 a b g r y_0 x_0^3&\\
    &+&8 a b r y_0 \alpha  x_0^3-12 b^2 g^2 r^2 \beta  x_0^3+44 b g r^2 \beta  x_0^3-12 r^2 \beta  x_0^3-12 a b r y_0 \beta  x_0^3&\\
    &+&4 a y_0 \mu  x_0^3-14 b g r \alpha  \mu  x_0^3+6 r \alpha  \mu  x_0^3+10 b g r \beta  \mu  x_0^3-4 r \beta  \mu  x_0^3
    &\\
    &+&6 b r \alpha  \beta  \mu  x_0^3+4 g^2 r^2 x_0^2+4 a^2 y_0^2 x_0^2-8 a c \alpha ^2 x_0^2+4 b^2 g^2 r^2 \beta ^2 x_0^2-14 b g r^2 \beta ^2 x_0^2
    &\\
    &+&4 r^2 \beta ^2 x_0^2+g^2 \mu ^2 x_0^2+\alpha ^2 \mu ^2 x_0^2+4 g \alpha  \mu ^2 x_0^2-8 a g r y_0 x_0^2-12 a c g \alpha  x_0^2&\\
    &-&6 a r y_0 \alpha  x_0^2+6 a b g r y_0 \alpha  x_0^2+14 b g^2 r^2 \beta  x_0^2-14 g r^2 \beta  x_0^2
    &\\
    &+&8 a r y_0 \beta  x_0^2-8 a b g r y_0 \beta  x_0^2-6 a b r y_0 \alpha  \beta  x_0^2+4 g^2 r \mu  x_0^2
    &\\
    &-&6 b g^2 r \alpha  \mu  x_0^2+10 g r \alpha  \mu  x_0^2+6 a y_0 \alpha  \mu  x_0^2+4 b g^2 r \beta  \mu  x_0^2&\\
    &+&6 g r \beta  \mu  x_0^2+10 b g r \alpha  \beta  \mu  x_0^2-4 r \alpha  \beta  \mu  x_0^2-4 a c g \alpha ^2 x_0&\\
    &-&4 b g^2 r^2 \beta ^2 x_0+4 g r^2 \beta ^2 x_0+2 g \alpha ^2 \mu ^2 x_0+2 g^2 \alpha  \mu ^2 x_0&\\
    &-&4 a g r y_0 \alpha  x_0-4 g^2 r^2 \beta  x_0+4 a g r y_0 \beta  x_0+4 a r y_0 \alpha  \beta  x_0&\\
    &+&2 a y_0 \alpha ^2 \mu  x_0+4 g^2 r \alpha  \mu  x_0-2 a g y_0 \alpha  \mu  x_0-2 g^2 r \beta  \mu  x_0&\\
    &+&4 b g^2 r \alpha  \beta  \mu  x_0-6 g r \alpha  \beta  \mu  x_0+a^2 y_0^2 \alpha ^2+g^2 r^2 \beta ^2&\\
    &+&g^2 \alpha ^2 \mu ^2+2 a g r y_0 \alpha  \beta -2 a g y_0 \alpha ^2 \mu -2 g^2 r \alpha  \beta  \mu&\\
    &+&10 g r \mu  x_0^3-34 b g r^2 x_0^4-14 b g r \mu  x_0^4-20 a c \alpha  x_0^3-34 b^2 g r^2 \beta  x_0^4&\\
    &-&8 b r \mu  x_0^5+9 b^2 g^2 r^2 x_0^4-6 b g^2 r \mu  x_0^3-4 a b g r y_0 \alpha  \beta  x_0&\\

\end{array}
\end{equation}
The following equations are derived from the critical point $P_2$ given by expression (\ref{PC1})  and are as follows.
\begin{equation}\label{ETA}
\begin{array}{lcll}
    \eta&=&2 \sqrt[3]{2} a^2 c^2+4 \sqrt[3]{2} a b c g \mu  r&\\
    &-&6 \sqrt[3]{2} a \alpha  b c \mu  r-4 \sqrt[3]{2} a b \beta  c \mu  r-4 \sqrt[3]{2} a c \mu  r+2 \sqrt[3]{2} b^2 g^2 \mu ^2 r^2&\\
    &+&2 \sqrt[3]{2} b^2 \beta  g \mu ^2 r^2+2 \sqrt[3]{2} b^2 \beta ^2 \mu ^2 r^2+2 \sqrt[3]{2} b g \mu ^2 r^2&\\
    &-&2 b \Gamma  g \mu  r-2 \sqrt[3]{2} b \beta  \mu ^2 r^2+2 b \beta  \Gamma  \mu  r+2^{2/3} \Gamma ^2+2 \sqrt[3]{2} \mu ^2 r^2+2 \Gamma  \mu  r
\end{array}
\end{equation}

\begin{equation}\label{GAMMA}
\begin{array}{lcll}
    \Gamma^3&=&-2 a^3 c^3-6 a^2 b c^2 g \mu  r+9 a^2 \alpha  b c^2 \mu  r+6 a^2 b \beta  c^2 \mu  r+6 a^2 c^2 \mu  r&\\
    &-&6 a b^2 c g^2 \mu ^2 r^2+9 a \alpha  b^2 c g \mu ^2 r^2+3 a b^2 \beta  c g \mu ^2 r^2&\\
    &-&9 a \alpha  b^2 \beta  c \mu ^2 r^2-6 a b^2 \beta ^2 c \mu ^2 r^2+3 a b c g \mu ^2 r^2-9 a \alpha  b c \mu ^2 r^2&\\
    &-&3 a b \beta  c \mu ^2 r^2-6 a c \mu ^2 r^2-2 b^3 g^3 \mu ^3 r^3-3 b^3 \beta  g^2 \mu ^3 r^3&\\
    &+&3 b^3 \beta ^2 g \mu ^3 r^3+2 b^3 \beta ^3 \mu ^3 r^3-3 b^2 g^2 \mu ^3 r^3-12 b^2 \beta  g \mu^3 r^3&\\
    &-&3 b^2 \beta ^2 \mu ^3 r^3+3 b g \mu ^3 r^3-3 b \beta  \mu ^3 r^3+2 \mu ^3 r^3 +\sqrt{R_{\Gamma}},
\end{array}
\end{equation}
where
\begin{equation}\label{RGAMMA}
\begin{array}{lcll}
    R_{\Gamma}&=&-27 b^4 g^4 r^6 \mu ^6-54 b^3 g^3 r^6 \mu ^6-27 b^2 g^2 r^6 \mu ^6&\\
    &-&27 b^4 r^6 \beta ^4 \mu ^6-27 b^6 g^2 r^6 \beta ^4 \mu ^6-54 b^5 g r^6 \beta ^4 \mu ^6+54 b^3 r^6 \beta ^3 \mu ^6-54 b^6 g^3 r^6 \beta ^3 \mu ^6&\\
    &-&54 b^5 g^2 r^6 \beta ^3 \mu ^6+54 b^4 g r^6 \beta ^3 \mu ^6-27 b^6 g^4 r^6 \beta^2 \mu^6&\\
    &+&54 b^5 g^3 r^6 \beta ^2 \mu ^6-27 b^2 r^6 \beta ^2 \mu ^6+162 b^4 g^2 r^6 \beta^2 \mu^6+54 b^3 g r^6 \beta^2 \mu^6&\\
    &+&54 b^5 g^4 r^6 \beta  \mu ^6+54 b^4 g^3 r^6 \beta  \mu ^6-54 b^3 g^2 r^6 \beta  \mu ^6-54 b^2 g r^6 \beta  \mu ^6&\\
    &-&54 a b^3 c g^3 r^5 \mu ^5+54 a b^2 c g^2 r^5 \mu ^5+54 a b^5 c g^2 r^5 \beta ^3 \mu ^5+54 a b^3 c r^5 \beta ^3 \mu^5&\\
    &+&216 a b^4 c g r^5 \beta ^3 \mu ^5-54 a b^4 c r^5 \alpha  \beta ^3 \mu ^5+54 a b^5 c g r^5 \alpha  \beta ^3 \mu^5&\\
    &-&54 a b^5 c g^3 r^5 \beta^2 \mu^5+108 a b^4 c g^2 r^5 \beta ^2 \mu ^5+54 a b^2 c r^5 \beta^2 \mu^5&\\
    &-&108 a b^3 c g r^5 \beta ^2 \mu ^5+216 a b^5 c g^2 r^5 \alpha  \beta ^2 \mu ^5+216 a b^3 c r^5 \alpha  \beta^2 \mu^5&\\
    &+&108 a b^4 c g r^5 \alpha  \beta ^2 \mu ^5+54 a b^4 c g^3 r^5 \alpha  \mu ^5+216 a b^3 c g^2 r^5 \alpha  \mu^5&\\
    &+&54 a b^2 c g r^5 \alpha  \mu ^5+216 a b^4 c g^3 r^5 \beta  \mu ^5+108 a b^3 c g^2 r^5 \beta  \mu^5&\\
    &+&216 a b^2 c g r^5 \beta  \mu^5+54 a b^5 c g^3 r^5 \alpha  \beta  \mu ^5-108 a b^4 c g^2 r^5 \alpha  \beta  \mu ^5-54 a b^2 c r^5 \alpha  \beta  \mu^5&\\
    &+&108 a b^3 c g r^5 \alpha  \beta  \mu ^5-27 a^2 b^2 c^2 g^2 r^4 \mu ^4-27 a^2 b^2 c^2 r^4 \alpha^2 \mu^4&\\
    &-&27 a^2 b^4 c^2 g^2 r^4 \alpha ^2 \mu ^4-270 a^2 b^3 c^2 g r^4 \alpha ^2 \mu ^4-27 a^2 b^2 c^2 r^4 \beta^2 \mu^4&\\
    &-&27 a^2 b^4 c^2 g^2 r^4 \beta ^2 \mu ^4-270 a^2 b^3 c^2 g r^4 \beta^2 \mu ^4-27 a^2 b^4 c^2 r^4 \alpha^2 \beta^2 \mu^4&\\
    &+&108 a^2 b^3 c^2 r^4 \alpha  \beta^2 \mu^4-108 a^2 b^4 c^2 g r^4 \alpha  \beta^2 \mu^4&\\
    &+&108 a^2 b^3 c^2 g^2 r^4 \alpha  \mu^4-108 a^2 b^2 c^2 g r^4 \alpha  \mu ^4+270 a^2 b^3 c^2 g^2 r^4 \beta  \mu^4&\\
    &-&270 a^2 b^2 c^2 g r^4 \beta  \mu ^4+270 a^2 b^3 c^2 r^4 \alpha^2 \beta  \mu^4-270 a^2 b^4 c^2 g r^4 \alpha^2 \beta  \mu^4&\\
    &+&108 a^2 b^2 c^2 r^4 \alpha  \beta  \mu^4+108 a^2 b^4 c^2 g^2 r^4 \alpha  \beta  \mu^4&\\
    &-&810 a^2 b^3 c^2 g r^4 \alpha  \beta  \mu^4+108 a^3 b^3 c^3 r^3 \alpha^3 \mu^3+54 a^3 b^2 c^3 r^3 \alpha ^2 \mu^3&\\
    &-&54 a^3 b^3 c^3 g r^3 \alpha^2 \mu^3+54 a^3 b^2 c^3 g r^3 \alpha  \mu^3+108 a^3 b^2 c^3 g r^3 \beta  \mu^3&\\
    &+&54 a^3 b^3 c^3 r^3 \alpha^2 \beta  \mu^3-54 a^3 b^2 c^3 r^3 \alpha  \beta  \mu^3+54 a^3 b^3 c^3 g r^3 \alpha  \beta  \mu^3-27 a^4 b^2 c^4 r^2 \alpha^2 \mu^2
\end{array}
\end{equation}

Here is the equation that we obtained to demonstrate's process in Proposition \ref{propositionSN1}:
\begin{equation}\label{S_N_1}
\begin{array}{lcll}
        SN&=& a^4 \alpha^2 c^4-4 a^3 b r \alpha^3 \mu  c^3-2 a^3 r \alpha^2 \mu  c^3+2 a^3 b g r \alpha^2 \mu  c^3\\
        &-&2 a^3 g r \alpha  \mu  c^3-2 a^3 b r \alpha^2 \beta  \mu  c^3-4 a^3 g r \beta  \mu  c^3\\
        &+&2 a^3 r \alpha  \beta  \mu  c^3-2 a^3 b g r \alpha  \beta  \mu  c^3+a^2 g^2 r^2 \mu ^2 c^2\\
        &+&a^2 r^2 \alpha ^2 \mu ^2 c^2+a^2 b^2 g^2 r^2 \alpha ^2 \mu ^2 c^2+10 a^2 b g r^2 \alpha ^2 \mu ^2 c^2\\
        &+&a^2 r^2 \beta ^2 \mu ^2 c^2+a^2 b^2 g^2 r^2 \beta ^2 \mu ^2 c^2+10 a^2 b g r^2 \beta ^2 \mu ^2 c^2\\
        &+&a^2 b^2 r^2 \alpha ^2 \beta ^2 \mu ^2 c^2-4 a^2 b r^2 \alpha  \beta ^2 \mu ^2 c^2+4 a^2 b^2 g r^2 \alpha  \beta ^2 \mu ^2 c^2\\
        &-&4 a^2 b g^2 r^2 \alpha  \mu ^2 c^2+4 a^2 g r^2 \alpha  \mu ^2 c^2-10 a^2 b g^2 r^2 \beta  \mu ^2 c^2
        \\
        &+&10 a^2 g r^2 \beta  \mu ^2 c^2-10 a^2 b r^2 \alpha ^2 \beta  \mu ^2 c^2+10 a^2 b^2 g r^2 \alpha ^2 \beta  \mu ^2 c^2-4 a^2 r^2 \alpha  \beta  \mu ^2 c^2
        \\
        &-&4 a^2 b^2 g^2 r^2 \alpha  \beta  \mu ^2 c^2+30 a^2 b g r^2 \alpha  \beta  \mu ^2 c^2+2 a b g^3 r^3 \mu ^3 c\\
        &-&2 a g^2 r^3 \mu ^3 c-2 a b^3 g^2 r^3 \beta ^3 \mu ^3 c-2 a b r^3 \beta ^3 \mu ^3 c\\
        &-&8 a b^2 g r^3 \beta ^3 \mu ^3 c+2 a b^2 r^3 \alpha  \beta ^3 \mu ^3 c-2 a b^3 g r^3 \alpha  \beta ^3 \mu ^3 c+2 a b^3 g^3 r^3 \beta ^2 \mu ^3 c\\
        &-&4 a b^2 g^2 r^3 \beta ^2 \mu ^3 c-2 a r^3 \beta ^2 \mu ^3 c+4 a b g r^3 \beta ^2 \mu ^3 c\\
        &-&8 a b^3 g^2 r^3 \alpha  \beta ^2 \mu ^3 c-8 a b r^3 \alpha  \beta ^2 \mu ^3 c-4 a b^2 g r^3 \alpha  \beta ^2 \mu ^3 c\\
        &-&2 a b^2 g^3 r^3 \alpha  \mu ^3 c-8 a b g^2 r^3 \alpha  \mu ^3 c-2 a g r^3 \alpha  \mu ^3 c\\
        &-&8 a b^2 g^3 r^3 \beta  \mu ^3 c-4 a b g^2 r^3 \beta  \mu ^3 c-8 a g r^3 \beta  \mu ^3 c\\
        &-&2 a b^3 g^3 r^3 \alpha  \beta  \mu ^3 c+4 a b^2 g^2 r^3 \alpha  \beta  \mu ^3 c+2 a r^3 \alpha  \beta  \mu ^3 c-4 a b g r^3 \alpha  \beta  \mu ^3 c\\
        &+&b^2 g^4 r^4 \mu ^4+2 b g^3 r^4 \mu ^4+g^2 r^4 \mu ^4+b^2 r^4 \beta ^4 \mu ^4+b^4 g^2 r^4 \beta ^4 \mu ^4\\
        &+&2 b^3 g r^4 \beta ^4 \mu ^4+2 b^4 g^3 r^4 \beta ^3 \mu ^4+2 b^3 g^2 r^4 \beta ^3 \mu ^4\\
        &-&2 b r^4 \beta ^3 \mu ^4-2 b^2 g r^4 \beta ^3 \mu ^4+b^4 g^4 r^4 \beta ^2 \mu ^4\\
        &-&2 b^3 g^3 r^4 \beta ^2 \mu ^4-6 b^2 g^2 r^4 \beta ^2 \mu^4-2 b g r^4 \beta ^2 \mu ^4+r^4 \beta ^2 \mu^4\\
        &-&2 b^3 g^4 r^4 \beta  \mu^4-2 b^2 g^3 r^4 \beta  \mu ^4+2 b g^2 r^4 \beta  \mu ^4+2 g r^4 \beta  \mu ^4,
\end{array}
\end{equation}

We derived the following expression in proof of Theorem \ref{BTToe}.

\begin{equation}\label{BT_1}
\begin{array}{lcll}
    BT_2&=&a^2 b^2 \beta^3 c^2 r^3-2 a^2 b \beta ^2 c^2 r^3-4 a^2 \alpha  b \beta  c^2 \mu  r^2-2 a^2 b \beta ^2 c^2 \mu  r^2\\
    &-&a^2 \alpha ^2 b c^2 \mu ^2 r-a^2 \alpha  b \beta  c^2 \mu ^2 r+a^2 \beta  c^2 r^3-2 a^2 \beta  c^2 \mu  r^2\\
    &-&a^2 \alpha  c^2 \mu ^2 r+a \alpha  b^2 \beta ^2 c \mu ^3 r^2-2 a \alpha  b \beta  c \mu ^3 r^2-2 a \alpha ^2 b c \mu ^4 r\\
    &+&a \alpha  c \mu ^5+a \alpha  c \mu ^3 r^2+2 a \beta  c \mu ^4 r+\alpha ^2 (-b) \mu ^6 r-\alpha  b \beta  \mu ^6 r-\alpha  \mu ^6 r-\beta  \mu ^6 r\\
\end{array}
\end{equation}

As a previous expression, the following equation gets provided in the proof of Theorem \ref{HopfTeo}.
\begin{equation}\label{H_1}
\rotatebox{270}{$\begin{array}{lcll}&\\
    HOPF&=& -b^2 \beta ^5 \mu ^3 r^5-b^4 g^2 \beta ^5 \mu ^3 r^5-2 b^3 g \beta ^5 \mu ^3 r^5-2 b^4 g^3 \beta ^4 \mu ^3 r^5-2 b^3 g^2 \beta ^4 \mu ^3 r^5+2 b \beta ^4 \mu ^3 r^5+2 b^2 g \beta ^4 \mu ^3 r^5-b^4 g^4 \beta ^3 \mu ^3 r^5&\\
    &+&6 b^2 g^2 \beta ^3 \mu ^3 r^5+2 b g \beta ^3 \mu ^3 r^5-\beta ^3 \mu ^3 r^5+2 b^3 g^4 \beta ^2 \mu ^3 r^5+2 b^2 g^3 \beta ^2 \mu ^3 r^5-2 b g^2 \beta ^2 \mu ^3 r^5-2 g \beta ^2 \mu ^3 r^5-b^2 g^4 \beta  \mu ^3 r^5&\\
    &-&g^2 \beta  \mu ^3 r^5+b^2 \alpha  \beta ^4 \mu ^4 r^4+b^4 g^2 \alpha  \beta ^4 \mu ^4 r^4+2 b^3 g \alpha  \beta ^4 \mu ^4 r^4+2 b^4 g^3 \alpha  \beta ^3 \mu ^4 r^4+2 b^3 g^2 \alpha  \beta ^3 \mu ^4 r^4-2 b \alpha  \beta ^3 \mu ^4 r^4&\\
    &-&2 b^2 g \alpha  \beta ^3 \mu ^4 r^4+b^4 g^4 \alpha  \beta ^2 \mu ^4 r^4-2 b^3 g^3 \alpha  \beta ^2 \mu ^4 r^4-6 b^2 g^2 \alpha  \beta ^2 \mu ^4 r^4-2 b g \alpha  \beta ^2 \mu ^4 r^4+\alpha  \beta ^2 \mu ^4 r^4-b^2 g^4 \alpha  \mu ^4 r^4&\\
    &+&2 b g^3 \alpha  \mu ^4 r^4+g^2 \alpha  \mu ^4 r^4-2 b^3 g^4 \alpha  \beta  \mu ^4 r^4-2 b^2 g^3 \alpha  \beta  \mu ^4 r^4+2 b g^2 \alpha  \beta  \mu ^4 r^4+2 g \alpha  \beta  \mu ^4 r^4+2 a b^3 c g^2 \beta ^4 \mu ^2 r^4+2 a b c \beta ^4 \mu ^2 r^4&\\
    &+&8 a b^2 c g \beta ^4 \mu ^2 r^4-2 a b^2 c \alpha  \beta ^4 \mu ^2 r^4+2 a b^3 c g \alpha  \beta ^4 \mu ^2 r^4-2 a b^3 c g^3 \beta ^3 \mu ^2 r^4+4 a b^2 c g^2 \beta ^3 \mu ^2 r^4+2 a c \beta ^3 \mu ^2 r^4-4 a b c g \beta ^3 \mu ^2 r^4&\\
    &+&8 a b^3 c g^2 \alpha  \beta ^3 \mu ^2 r^4+8 a b c \alpha  \beta ^3 \mu ^2 r^4+4 a b^2 c g \alpha  \beta ^3 \mu ^2 r^4+8 a b^2 c g^3 \beta ^2 \mu ^2 r^4+4 a b c g^2 \beta ^2 \mu ^2 r^4+8 a c g \beta ^2 \mu ^2 r^4+2 a b^3 c g^3 \alpha  \beta ^2 \mu ^2 r^4&\\
    &-&4 a b^2 c g^2 \alpha  \beta ^2 \mu ^2 r^4-2 a c \alpha  \beta ^2 \mu ^2 r^4+4 a b c g \alpha  \beta ^2 \mu ^2 r^4-2 a b c g^3 \beta  \mu ^2 r^4+2 a c g^2 \beta  \mu ^2 r^4+2 a b^2 c g^3 \alpha  \beta  \mu ^2 r^4+8 a b c g^2 \alpha  \beta  \mu ^2 r^4&\\
    &+&2 a c g \alpha  \beta  \mu ^2 r^4+b^2 g^2 \beta ^3 \mu ^5 r^3+2 b g \beta ^3 \mu ^5 r^3+b^3 g^2 \alpha  \beta ^3 \mu ^5 r^3+b \alpha  \beta ^3 \mu ^5 r^3+2 b^2 g \alpha  \beta ^3 \mu ^5 r^3+\beta ^3 \mu ^5 r^3+b^3 g^4 \alpha ^2 \mu ^5 r^3&\\
    &+&2 b^2 g^3 \alpha ^2 \mu ^5 r^3+b g^2 \alpha ^2 \mu ^5 r^3+2 b^2 g^3 \beta ^2 \mu ^5 r^3+4 b g^2 \beta ^2 \mu ^5 r^3+b^3 g^2 \alpha ^2 \beta ^2 \mu ^5 r^3+b \alpha ^2 \beta ^2 \mu ^5 r^3+2 b^2 g \alpha ^2 \beta ^2 \mu ^5 r^3+2 g \beta ^2 \mu ^5 r^3&\\
    &+&2 b^3 g^3 \alpha  \beta ^2 \mu ^5 r^3+5 b^2 g^2 \alpha  \beta ^2 \mu ^5 r^3+4 b g \alpha  \beta ^2 \mu ^5 r^3+\alpha  \beta ^2 \mu ^5 r^3+b^2 g^4 \alpha  \mu ^5 r^3+2 b g^3 \alpha  \mu ^5 r^3+g^2 \alpha  \mu ^5 r^3+b^2 g^4 \beta  \mu ^5 r^3&\\
    &+&2 b g^3 \beta  \mu ^5 r^3+g^2 \beta  \mu ^5 r^3+2 b^3 g^3 \alpha ^2 \beta  \mu ^5 r^3+4 b^2 g^2 \alpha ^2 \beta  \mu ^5 r^3+2 b g \alpha ^2 \beta  \mu ^5 r^3+b^3 g^4 \alpha  \beta  \mu ^5 r^3+4 b^2 g^3 \alpha  \beta  \mu ^5 r^3+5 b g^2 \alpha  \beta  \mu ^5 r^3&\\
    &+&2 g \alpha  \beta  \mu ^5 r^3-a b^2 c g^2 \beta ^3 \mu ^3 r^3+a b^2 c \alpha ^2 \beta ^3 \mu ^3 r^3-3 a b^3 c g \alpha ^2 \beta ^3 \mu ^3 r^3-2 a c \beta ^3 \mu ^3 r^3-3 a b c g \beta ^3 \mu ^3 r^3-2 a b^3 c g^2 \alpha  \beta ^3 \mu ^3 r^3&\\
    &&4 a b c \alpha  \beta ^3 \mu ^3 r^3-10 a b^2 c g \alpha  \beta ^3 \mu ^3 r^3-3 a b^2 c g^3 \alpha ^2 \mu ^3 r^3-10 a b c g^2 \alpha ^2 \mu ^3 r^3-3 a c g \alpha ^2 \mu ^3 r^3-a b^2 c g^3 \beta ^2 \mu ^3 r^3-4 a b c g^2 \beta ^2 \mu ^3 r^3&\\
    &-&10 a b^3 c g^2 \alpha ^2 \beta ^2 \mu ^3 r^3-8 a b c \alpha ^2 \beta ^2 \mu ^3 r^3-6 a b^2 c g \alpha ^2 \beta ^2 \mu ^3 r^3-3 a c g \beta ^2 \mu ^3 r^3+2 a b^3 c g^3 \alpha  \beta ^2 \mu ^3 r^3-10 a b^2 c g^2 \alpha  \beta ^2 \mu ^3 r^3&\\
    &-&4 a c \alpha  \beta ^2 \mu ^3 r^3-4 a b c g \alpha  \beta ^2 \mu ^3 r^3+2 a b c g^3 \alpha  \mu ^3 r^3-2 a c g^2 \alpha  \mu ^3 r^3-a b c g^3 \beta  \mu ^3 r^3-a c g^2 \beta  \mu ^3 r^3-3 a b^3 c g^3 \alpha ^2 \beta  \mu ^3 r^3&\\
    &+&2 a b^2 c g^2 \alpha ^2 \beta  \mu ^3 r^3+a c \alpha ^2 \beta  \mu ^3 r^3-6 a b c g \alpha ^2 \beta  \mu ^3 r^3-12 a b^2 c g^3 \alpha  \beta  \mu ^3 r^3-10 a b c g^2 \alpha  \beta  \mu ^3 r^3-10 a c g \alpha  \beta  \mu ^3 r^3-9 a^2 b c^2 g \beta ^3 \mu  r^3&\\
\end{array}$} 
\end{equation}   
    
\begin{equation}
\nonumber
\rotatebox{270}{$\begin{array}{lcll}&\\

    &+&5 a^2 b c^2 \alpha  \beta ^3 \mu  r^3-5 a^2 b^2 c^2 g \alpha  \beta ^3 \mu  r^3+9 a^2 b c^2 g^2 \beta ^2 \mu  r^3+9 a^2 b c^2 \alpha ^2 \beta ^2 \mu  r^3-9 a^2 b^2 c^2 g \alpha ^2 \beta ^2 \mu  r^3-9 a^2 c^2 g \beta ^2 \mu  r^3&\\
    &+&5 a^2 c^2 \alpha  \beta ^2 \mu  r^3+5 a^2 b^2 c^2 g^2 \alpha  \beta ^2 \mu  r^3-24 a^2 b c^2 g \alpha  \beta ^2 \mu  r^3-9 a^2 b c^2 g \alpha ^2 \beta  \mu  r^3+5 a^2 b c^2 g^2 \alpha  \beta  \mu  r^3-5 a^2 c^2 g \alpha  \beta  \mu  r^3&\\
    &-&6 a b^2 c g^2 \alpha ^3 \mu ^4 r^2-a c \alpha ^3 \mu ^4 r^2-3 a b c g \alpha ^3 \mu ^4 r^2+3 a b^2 c g^3 \alpha ^2 \mu ^4 r^2-6 a b c g^2 \alpha ^2 \mu ^4 r^2-a c g \alpha ^2 \mu ^4 r^2-a b^2 c \alpha ^3 \beta ^2 \mu ^4 r^2&\\
    &-&3 a c g \beta ^2 \mu ^4 r^2-a b^2 c g^2 \alpha  \beta ^2 \mu ^4 r^2-6 a b c g \alpha  \beta ^2 \mu ^4 r^2+3 a b c g^3 \alpha  \mu ^4 r^2-a c g^2 \alpha  \mu ^4 r^2+3 a b c g^3 \beta  \mu ^4 r^2+2 a b c \alpha ^3 \beta  \mu ^4 r^2&\\
    &-&6 a b^2 c g^2 \alpha ^2 \beta  \mu ^4 r^2-13 a b c g \alpha ^2 \beta  \mu ^4 r^2+3 a b^2 c g^3 \alpha  \beta  \mu ^4 r^2-4 a b c g^2 \alpha  \beta  \mu ^4 r^2-6 a c g \alpha  \beta  \mu ^4 r^2+2 a^3 c^3 g \beta ^2 r^2-2 a^3 c^3 \alpha  \beta ^2 r^2&\\
    &+&2 a^3 b c^3 g \alpha  \beta ^2 r^2+2 a^2 c^2 \alpha ^3 \mu ^2 r^2+2 a^2 b^2 c^2 g^2 \alpha ^3 \mu ^2 r^2+17 a^2 b c^2 g \alpha ^3 \mu ^2 r^2-7 a^2 b c^2 g^2 \alpha ^2 \mu ^2 r^2+7 a^2 c^2 g \alpha ^2 \mu ^2 r^2+2 a^2 b^2 c^2 \alpha ^3 \beta ^2 \mu ^2 r^2&\\
    &-&2 a^2 b c^2 g^2 \beta ^2 \mu ^2 r^2-2 a^2 b c^2 \alpha ^2 \beta ^2 \mu ^2 r^2+7 a^2 b^2 c^2 g \alpha ^2 \beta ^2 \mu ^2 r^2+5 a^2 c^2 g \beta ^2 \mu ^2 r^2+a^2 c^2 \alpha  \beta ^2 \mu ^2 r^2&\\
    &-&8 a^2 b c^2 \alpha ^3 \beta  \mu ^2 r^2+17 a^2 b^2 c^2 g \alpha ^3 \beta  \mu ^2 r^2-2 a^2 c^2 g^2 \beta  \mu ^2 r^2-2 a^2 c^2 \alpha ^2 \beta  \mu ^2 r^2-7 a^2 b^2 c^2 g^2 \alpha ^2 \beta  \mu ^2 r^2+49 a^2 b c^2 g \alpha ^2 \beta  \mu ^2 r^2&\\
    &+&18 a^2 c^2 g \alpha  \beta  \mu ^2 r^2+2 a^3 c^3 g \alpha  \beta  r^2-a b c \alpha ^4 \mu ^5 r-a c \alpha ^3 \mu ^5 r+2 a b c g \alpha ^3 \mu ^5 r-a b c g^2 \alpha ^2 \mu ^5 r+2 a c g \alpha ^2 \mu ^5 r-a c g^2 \alpha  \mu ^5 r&\\
    &-&a c \alpha ^2 \beta  \mu ^5 r+2 a b c g \alpha ^2 \beta  \mu ^5 r-a b c g^2 \alpha  \beta  \mu ^5 r+2 a c g \alpha  \beta  \mu ^5 r+6 a^2 b c^2 \alpha ^4 \mu ^3 r+4 a^2 c^2 \alpha ^3 \mu ^3 r-5 a^2 b c^2 g \alpha ^3 \mu ^3 r&\\
    &+&2 a^2 c^2 g^2 \alpha  \mu ^3 r+4 a^2 b c^2 \alpha ^3 \beta  \mu ^3 r+2 a^2 c^2 g^2 \beta  \mu ^3 r+2 a^2 c^2 \alpha ^2 \beta  \mu ^3 r-3 a^2 b c^2 g \alpha ^2 \beta  \mu ^3 r+2 a^2 b c^2 g^2 \alpha  \beta  \mu ^3 r-a^2 c^2 g \alpha  \beta  \mu ^3 r&\\
&+&2 b^3 g^3 \beta ^3 \mu ^3 r^5-2 b g^3 \beta  \mu ^3 r^5-6 a^3 c^3 g \alpha  \beta  \mu  r+2 a^2 b c^2 g^2 \alpha ^2 \mu ^3 r-3 a^2 c^2 g \alpha ^2 \mu ^3 r-a b c \alpha ^3 \beta  \mu ^5 r-a c g^2 \beta  \mu ^5 r&\\
    
    &-&9 a^3 b c^3 \alpha ^4 \mu  r-5 a^3 c^3 \alpha ^3 \mu  r+5 a^3 b c^3 g \alpha ^3 \mu  r-2 a^3 c^3 g \alpha ^2 \mu  r-5 a^3 b c^3 \alpha ^3 \beta  \mu  r-2 a^3 b c^3 g \alpha ^2 \beta  \mu  r+2 a^4 c^4 \alpha ^3-a^3 c^3 \alpha^3 \mu^2&\\
    &-&19 a^2 b c^2 g^2 \alpha  \beta  \mu ^2 r^2+18 a^2 b c^2 g \alpha  \beta ^2 \mu ^2 r^2-3 a b^2 c g \alpha ^3 \beta  \mu ^4 r^2-a b^2 c g \alpha ^2 \beta ^2 \mu ^4 r^2
\end{array}$}
\nonumber
\end{equation}

Next we present the complete equation for describing Generalized Hopf bifurcation.
\begin{equation}\label{Bau_1}
\rotatebox{270}{$\begin{array}{lcll}
    BAU&=&-24640 a^{14} b^3 g^2 r^3 \alpha ^{13} c^{14}-24640 a^{14} b^2 g r^3 \alpha ^{13} c^{14}+11200 a^{14} b^3 g^3 r^3 \alpha ^{12} c^{14}-81200 a^{14} b^2 g^2 r^3 \alpha ^{12} c^{14}-64680 a^{14} b g r^3 \alpha ^{12} c^{14}&\\
    &+&38920 a^{14} b^2 g^3 r^3 \alpha ^{11} c^{14}-57120 a^{14} b g^2 r^3 \alpha ^{11} c^{14}-40320 a^{14} g r^3 \alpha ^{11} c^{14}+28000 a^{14} b g^3 r^3 \alpha ^{10} c^{14}+28000 a^{14} b^3 g^3 r^3 \alpha ^{10} \beta ^2 c^{14}&\\
    &+&114240 a^{14} b^2 g^2 r^3 \alpha ^{10} \beta ^2 c^{14}-306600 a^{14} b g r^3 \alpha ^{10} \beta ^2 c^{14}+56000 a^{14} b^2 g^3 r^3 \alpha ^9 \beta ^2 c^{14}-57120 a^{14} b g^2 r^3 \alpha ^9 \beta ^2 c^{14}&\\
    &-&120960 a^{14} g r^3 \alpha ^9 \beta ^2 c^{14}+28000 a^{14} b g^3 r^3 \alpha ^8 \beta ^2 c^{14}-24640 a^{14} b^3 g r^3 \alpha ^{13} \beta  c^{14}-81200 a^{14} b^3 g^2 r^3 \alpha ^{12} \beta  c^{14}&\\
    &-&154000 a^{14} b^2 g r^3 \alpha ^{12} \beta  c^{14}+38920 a^{14} b^3 g^3 r^3 \alpha ^{11} \beta  c^{14}-195440 a^{14} b^2 g^2 r^3 \alpha ^{11} \beta  c^{14}-250320 a^{14} b g r^3 \alpha ^{11} \beta  c^{14}
    &\\
    &+&94920 a^{14} b^2 g^3 r^3 \alpha ^{10} \beta  c^{14}-114240 a^{14} b g^2 r^3 \alpha ^{10} \beta  c^{14}-120960 a^{14} g r^3 \alpha ^{10} \beta  c^{14}+56000 a^{14} b g^3 r^3 \alpha ^9 \beta  c^{14}
    &\\
    &-&225120 a^{13} b^3 g^3 r^5 \alpha ^{10} \beta ^2 c^{13}-455140 a^{13} b^2 g^2 r^5 \alpha ^{10} \beta ^2 c^{13}-161280 a^{13} b g r^5 \alpha ^{10} \beta ^2 c^{13}+133560 a^{13} b^3 g^4 r^5 \alpha ^9 \beta ^2 c^{13}
    &\\
    &-&455560 a^{13} b^2 g^3 r^5 \alpha ^9 \beta ^2 c^{13}-532980 a^{13} b g^2 r^5 \alpha ^9 \beta ^2 c^{13}-40320 a^{13} g r^5 \alpha ^9 \beta ^2 c^{13}+242760 a^{13} b^2 g^4 r^5 \alpha ^8 \beta ^2 c^{13}
    &\\
    &-&257040 a^{13} b g^3 r^5 \alpha ^8 \beta ^2 c^{13}-241920 a^{13} g^2 r^5 \alpha ^8 \beta ^2 c^{13}+126000 a^{13} b g^4 r^5 \alpha ^7 \beta ^2 c^{13}-24640 a^{13} b^3 g^2 r^5 \alpha ^{12} \beta  c^{13}
    &\\
    &+&24640 a^{13} b^2 g r^5 \alpha ^{12} \beta  c^{13}-25760 a^{13} b^3 g^3 r^5 \alpha ^{11} \beta  c^{13}-118160 a^{13} b^2 g^2 r^5 \alpha ^{11} \beta  c^{13}-64680 a^{13} b g r^5 \alpha ^{11} \beta  c^{13}
    &\\
    &+&16800 a^{13} b^3 g^4 r^5 \alpha ^{10} \beta  c^{13}-82880 a^{13} b^2 g^3 r^5 \alpha ^{10} \beta  c^{13}-154140 a^{13} b g^2 r^5 \alpha ^{10} \beta  c^{13}-40320 a^{13} g r^5 \alpha ^{10} \beta  c^{13}
    &\\
    &+&58380 a^{13} b^2 g^4 r^5 \alpha ^9 \beta  c^{13}-57680 a^{13} b g^3 r^5 \alpha ^9 \beta  c^{13}-60480 a^{13} g^2 r^5 \alpha ^9 \beta  c^{13}+42000 a^{13} b g^4 r^5 \alpha ^8 \beta  c^{13}+4480 a^{12} b^3 g^4 r^7 \alpha ^9 \beta ^2 c^{12}
    &\\
    &-&134120 a^{12} b^2 g^3 r^7 \alpha ^9 \beta ^2 c^{12}-97020 a^{12} b g^2 r^7 \alpha ^9 \beta ^2 c^{12}+5600 a^{12} b^3 g^5 r^7 \alpha ^8 \beta ^2 c^{12}+17780 a^{12} b^2 g^4 r^7 \alpha ^8 \beta ^2 c^{12}
    &\\
    &-&118020 a^{12} b g^3 r^7 \alpha ^8 \beta ^2 c^{12}-60480 a^{12} g^2 r^7 \alpha ^8 \beta ^2 c^{12}+19460 a^{12} b^2 g^5 r^7 \alpha ^7 \beta ^2 c^{12}+13440 a^{12} b g^4 r^7 \alpha ^7 \beta ^2 c^{12}
    &\\
    &\cdot& &\\
    &\cdot& &\\
    &\cdot& &\\
\end{array}$}
\end{equation}    

\begin{equation}
\nonumber
\rotatebox{270}{$\begin{array}{lcll}
    
    &-&144 b^6 g^8 r^{12} \beta ^{10} \mu ^{19}-672 b^5 g^7 r^{12} \beta ^{10} \mu ^{19}-1168 b^4 g^6 r^{12} \beta ^{10} \mu ^{19}-896 b^3 g^5 r^{12} \beta ^{10} \mu ^{19}-256 b^2 g^4 r^{12} \beta ^{10} \mu ^{19}
    &\\
    &-&16 b^{12} g^8 r^{12} \alpha ^6 \beta ^{10} \mu ^{19}-96 b^{11} g^7 r^{12} \alpha ^6 \beta ^{10} \mu ^{19}-208 b^{10} g^6 r^{12} \alpha ^6 \beta ^{10} \mu ^{19}-192 b^9 g^5 r^{12} \alpha ^6 \beta ^{10} \mu ^{19}
    &\\
    &-&64 b^8 g^4 r^{12} \alpha ^6 \beta ^{10} \mu ^{19}-160 b^{11} g^8 r^{12} \alpha ^5 \beta ^{10} \mu ^{19}-896 b^{10} g^7 r^{12} \alpha ^5 \beta ^{10} \mu ^{19}-1824 b^9 g^6 r^{12} \alpha ^5 \beta ^{10} \mu ^{19}
    &\\
    &-&1600 b^8 g^5 r^{12} \alpha ^5 \beta ^{10} \mu ^{19}-512 b^7 g^4 r^{12} \alpha ^5 \beta ^{10} \mu ^{19}-624 b^{10} g^8 r^{12} \alpha ^4 \beta ^{10} \mu ^{19}-3296 b^9 g^7 r^{12} \alpha ^4 \beta ^{10} \mu ^{19}
    &\\
    &-&6384 b^8 g^6 r^{12} \alpha ^4 \beta ^{10} \mu ^{19}-5376 b^7 g^5 r^{12} \alpha ^4 \beta ^{10} \mu ^{19}-1664 b^6 g^4 r^{12} \alpha ^4 \beta ^{10} \mu ^{19}-1216 b^9 g^8 r^{12} \alpha ^3 \beta ^{10} \mu ^{19}
    &\\
    &-&6144 b^8 g^7 r^{12} \alpha ^3 \beta ^{10} \mu ^{19}-11456 b^7 g^6 r^{12} \alpha ^3 \beta ^{10} \mu ^{19}-9344 b^6 g^5 r^{12} \alpha ^3 \beta ^{10} \mu ^{19}-2816 b^5 g^4 r^{12} \alpha ^3 \beta ^{10} \mu ^{19}
    &\\
    &-&1264 b^8 g^8 r^{12} \alpha ^2 \beta ^{10} \mu ^{19}-6176 b^7 g^7 r^{12} \alpha ^2 \beta ^{10} \mu ^{19}-11184 b^6 g^6 r^{12} \alpha ^2 \beta ^{10} \mu ^{19}-8896 b^5 g^5 r^{12} \alpha ^2 \beta ^{10} \mu ^{19}
    &\\
    &-&2624 b^4 g^4 r^{12} \alpha ^2 \beta ^{10} \mu ^{19}-672 b^7 g^8 r^{12} \alpha  \beta ^{10} \mu ^{19}-3200 b^6 g^7 r^{12} \alpha  \beta ^{10} \mu ^{19}-5664 b^5 g^6 r^{12} \alpha  \beta ^{10} \mu ^{19}
    &\\
    &-&4416 b^4 g^5 r^{12} \alpha  \beta ^{10} \mu ^{19}-1280 b^3 g^4 r^{12} \alpha  \beta ^{10} \mu ^{19}-80 b^{12} g^8 r^{12} \alpha ^7 \beta ^9 \mu ^{19}-480 b^{11} g^7 r^{12} \alpha ^7 \beta ^9 \mu ^{19}
    &\\
    &-&1040 b^{10} g^6 r^{12} \alpha ^7 \beta ^9 \mu ^{19}-960 b^9 g^5 r^{12} \alpha ^7 \beta ^9 \mu ^{19}-320 b^8 g^4 r^{12} \alpha ^7 \beta ^9 \mu ^{19}-56 b^{12} g^9 r^{12} \alpha ^6 \beta ^9 \mu ^{19}
    &\\
    &-&136 b^{11} g^8 r^{12} \alpha ^6 \beta ^9 \mu ^{19}-5208 b^{10} g^7 r^{12} \alpha ^6 \beta ^9 \mu ^{19}-9792 b^9 g^6 r^{12} \alpha ^6 \beta ^9 \mu ^{19}-8224 b^8 g^5 r^{12} \alpha ^6 \beta ^9 \mu ^{19}
    &\\
    &-&2560 b^7 g^4 r^{12} \alpha ^6 \beta ^9 \mu ^{19}-560 b^{11} g^9 r^{12} \alpha ^5 \beta ^9 \mu ^{19}-6256 b^{10} g^8 r^{12} \alpha ^5 \beta ^9 \mu ^{19}-22864 b^9 g^7 r^{12} \alpha ^5 \beta ^9 \mu ^{19}
    &\\
    &-&37520 b^8 g^6 r^{12} \alpha ^5 \beta ^9 \mu ^{19}+28672 b^7 g^5 r^{12} \alpha ^5 \beta ^9 \mu ^{19}
\end{array}$}
\nonumber
\end{equation}

\newpage
\bibliographystyle{unsrt}

\bibliography{mybibfile}

\begin{thebibliography}{10}

\bibitem{Altrock2015}
P.~L. Altrock and F.~Michor.
\newblock {The mathematics of cancer: Integrating quantitative models}.
\newblock {\em Nat. Rev. Cancer}, 15:730--745, 2015.

\bibitem{Kolev}
M.~Kolev.
\newblock {Mathematical modelling of the competition between tumors and immune
  system considering the role of the antibodies}.
\newblock {\em Mathematical and Computer Modelling}, 37(11):1143--1152, 2003.

\bibitem{Gore2009}
H.~Youk J.~Gore and A.~van Oudenaarden.
\newblock {Snowdrift game dynamics and facultative cheating in yeast}.
\newblock {\em Nature}, 459:253--256, 2009.

\bibitem{Axelrod2015}
D.~E.~Axelrod R.~Axelrod and K.~J. Pienta.
\newblock {Evolution of cooperation among tumor cells}.
\newblock {\em Nature}, 103:13474--13479, 2015.

\bibitem{Courchamp2008}
L.~Berec F.~Courchamp and J.~Gascoigne.
\newblock {\em {Allee effects in ecology and conservation}}.
\newblock OUP Oxford, 2008.

\bibitem{Bellomo}
N.~Bellomo and L.~Preziosi.
\newblock Modelling and mathematical problems related to tumor evolution and
  its interaction with the immune system.
\newblock {\em Mathematical and Computer Modelling}, 32(3):413--452, 2000.

\bibitem{Boer1985}
R.~J.~De Boer.
\newblock {Macrophage {T} lymphocyte interactions in the anti-tumor immune
  response: {A} mathematical model}.
\newblock {\em J. Immunol.}, 134:2748--2758, 1985.

\bibitem{Kuznetsov1994}
M.~A.~Taylor V.~A.~Kuznetsov, I. A.~Makalkin and A.~S. Perelson.
\newblock {Nonlinear dynamics of immunogenic tumors: Parameter estimation and
  global bifurcation analysis}.
\newblock {\em Bull. Math. Biol.}, 56(1):1--15, 1994.

\bibitem{Delgado2020}
E.~Hern\'andez-L\'opez J.~Delgado and L.~I. Hern\'andez-Mart\'inez.
\newblock {Bautin bifurcation in a minimal model of immunoediting}.
\newblock {\em Discr. Contin. Dyn. Syst. B}, 25:1397--1414, 2020.

\bibitem{Kirschner1998}
D.~Kirschner and J.~C. Panetta.
\newblock {Modeling immunotherapy of the tumor-immune interaction}.
\newblock {\em J. Math. Biol.}, 37:235--252, 1998.

\bibitem{Starkov2013}
K.~E. Starkov and L.~N. Coria.
\newblock {Global dynamics of the Kirschner \& Panetta model for the tumor
  immunotherapy}.
\newblock {\em Nonlin. Anal.: Real World Appl.}, 14:1425--1433., 2013.

\bibitem{Xavier2017}
R.~W. dos~Santos M.~P.~Xavier, C. R.~Bonin and M.~Lobosco.
\newblock {On the use of {G}illespie stochastic simulation algorithm in a model
  of the human immune system response to the yellow fever vaccine}.
\newblock {\em IEEE Int. Conf. Bioinformatics and Biomedicine (BIBM)}, pages
  1476--1482, 2017.

\bibitem{Wei2013}
H.~C. Wei and J.~T. Lin.
\newblock {Periodically pulsed immunotherapy in a mathematical model of tumor
  immune interaction}.
\newblock {\em Int. J. Bifurcation and Chaos}, pages 1350068:1--13, 2013.

\bibitem{hernandez2021bifurcations}
E.~Hern\'andez-L\'opez and M.~N\'unez-L\'opez.
\newblock {Bifurcations in a cancer and immune model with {A}llee effect}.
\newblock {\em International Journal of Bifurcation and Chaos},
  31(13):2130039:1--14, 2021.

\bibitem{Rocha2019}
J.~L. Rocha and A.~K. Taha.
\newblock {Allee effect bifurcation in generalized logistic maps}.
\newblock {\em Int. J. Bifurcation and Chaos}, 29:1950039:1--19, 2019.

\bibitem{Capistran2018}
M.~N\'u\~nez-L\'opez M.~A.~Capistr\'an and G.~A. Rempala.
\newblock {Extracellular dynamics of early HIV infection}.
\newblock {\em Math. Meth. Appl. Sci.}, 41:8859--8870, 2018.

\bibitem{Celik2009}
C.~Celik and O.~Duman.
\newblock {Allee effect in a discrete time predator--prey system}.
\newblock {\em Chaos Solit. Fract.}, pages 1956--1962, 2009.

\end{thebibliography}

\end{document}